\documentclass{article}
\usepackage{psfrag}
\usepackage{graphicx,amsfonts,amsmath}
\usepackage{pgfplots}
\usepackage[margin=1in]{geometry}
\usepackage[english]{babel}
\usepackage{csquotes}
\usepackage{authblk}
\usepackage{amsrefs}
\usepackage{xcolor}
\usepackage{float}
\usepackage{bigints}
\usepackage{tikz}
\usepackage{colonequals}
\newcommand\myeq{\mathrel{\overset{\makebox[0pt]{\mbox{\normalfont\tiny\sffamily \text{by} (\ref{Rankine2})}}}{=}}}
\numberwithin{equation}{section}
\usepackage[para]{footmisc}
\usetikzlibrary{intersections}
\usetikzlibrary{external}

\usepgfplotslibrary{fillbetween}
\pgfplotsset{width=10cm,compat=1.9}
\usepgfplotslibrary{external}
  
\title{An Analysis of a $2 \times 2$ Keyfitz-Kranzer Type Balance System with Varying Generalized Chaplygin Gas}
\date{\today}
\bibliography{plain}

\author[1]{Jack Frew \thanks{frew.25@osu.edu}}
\author[1]{Nigel Keyser \thanks{keyser.88@osu.edu}}
\author[2]{Ethan Kim \thanks{etk2124@columbia.edu}}
\author[3]{Griffin Paddock \thanks{griffinmilo@vt.edu}}
\author[4]{Camden Toumbleston \thanks{cstoumbl@ncsu.edu}}
\author[5]{Sara Wilson \thanks{srw81@pitt.edu}}
\author[6]{Charis Tsikkou \thanks{tsikkou@math.wvu.edu}}
\affil[1]{The Ohio State University, Columbus, OH, 43210, USA}
\affil[2]{School of Engineering and Applied Sciences, Columbia University, New York City, NY, 10027, USA}
\affil[3]{Virginia Tech, Blacksburg, Virginia, 24060, USA}
\affil[4]{North Carolina State University, Raleigh, NC, 27606, USA}
\affil[5]{University of Pittsburgh, Pittsburgh, PA, 15260, USA}
\affil[6]{School of Mathematical and Data Sciences, West Virginia University, Morgantown, WV 26506, USA}

\begin{document}

\maketitle
\begin{abstract}\noindent
We consider a system of two balance laws of Keyfitz-Kranzer type with varying generalized Chaplygin gas, which exhibits negative pressure and is a product of a function of time and the inverse of a power of the density. The Chaplygin gas is a fluid designed to accommodate measurements for the early universe and late-time universal expansion while obeying the pressure-density-time relation. We produce an explanation and description of the non-self-similar Riemann solutions, including the non-classical singular solutions. We also find that due to a direct dependence on time, a change in the regions allowing for combinations of classical and non-classical singular solutions occurs, therefore a Riemann solution can have different solutions over several time intervals. Our findings are confirmed numerically using the local Lax-Friederichs scheme. 
\end{abstract}

\vspace{5mm}


\vspace{5mm}


\section{Introduction}

It is well accepted that dark energy plays a critical role in the current expansion of the universe, in particular, the shift from deceleration to acceleration in the current epoch.  Connections between dark energy and dark matter are still debated in the physics community; the Chaplygin gas models are possible candidates for such connections, exhibiting early behavior akin to dark matter and later behavior akin to a cosmological constant.  All Chaplygin gas models are unified in describing dark energy through scalar fields with a negative pressure and inverse relation density of the form $p=A/\rho$, where $A<0$.  Of particular interest are the varying and generalized Chaplygin gas models \cite{Khurshudyan,Lipscombe}.  The former is characterized by the equation of state $p=B(t)/\rho$.  Recently, Khurshudyan \cite{Khurshudyan} proposed the idea of setting $B(t) = Ae^{\eta t}$, where $\eta$ and $A<0$ are constants. He showed that this form describes the quasi-exponential phase of the universe.
Li \cite{Li}, built upon this work by analyzing the system of balance laws formed from a mathematical viewpoint. The generalized Chaplygin gas is characterized by the equation of state $p=A\rho ^\gamma$, with $\gamma < 0$, as discussed in \cite{Sen}.  To the best of our knowledge, there has been little research on a combination of both the varying and generalized model, denoted by the varying generalized Chaplygin gas model (VGCG).

\vspace{5mm}

\noindent
This paper aims to study the solutions to the Riemann problem, an initial value problem that consists of data containing two constant states separated by a discontinuity at the origin, to a non-symmetric Keyfitz-Kranzer type system:
\begin{equation}\label{EQMAIN}
    \begin{split}\begin{cases}
    \rho_t+\bigg(\rho \bigg(u-p\big(\rho\big)\bigg)\bigg)_x&=k\rho,\\
    \big(\rho u\big)_t+\bigg(\rho u \bigg(u-p\big(\rho\big)\bigg)\bigg)_x&=\eta \rho u+\beta \rho, 
   \end{cases}
    \end{split}
\end{equation}
where $A<0,$ $k,$ $\eta,$ and $\beta$ are non-zero physical constants. The independent variables are time $t \in \mathbb{R}^+$ and position $x \in \mathbb{R},$ and the dependent variables are density $\rho$, and fluid velocity $u.$ We restrict attention to $p(\rho)=A \rho^{\gamma} e^{\eta t},$ with $\gamma<0,$ ($\gamma\neq -1$) and $\rho(x,t)>0.$ We refer the reader to Li \cite{Li} for the special case $\gamma= -1.$ 
In particular, we focus on the existence of singular solutions (delta or singular shocks) which denote the mass' concentration process and may be interpreted as galaxies in the universe. The non-autonomous system of balance laws (\ref{EQMAIN}) is of great interest mathematically and physically, as the solutions are non-self-similar and the shock and rarefaction curves change over time. To the best of our knowledge, direct time dependence resulting in changes in the areas where classical and non-classical singular solutions exist has not been analyzed and confirmed numerically before. 

\vspace{5mm}

The singular solutions involve the so-called delta or singular shocks, a more compressive generalization of the ordinary shock wave, where at least one state variable develops an extreme concentration in the form of a weighted Dirac delta function. They were initially discovered by Keyfitz and Kranzer \cite{Ke_Kr_1, Ke_Kr_3, Ke_Kr_2} and later studied in greater depth by Sever \cite{Se}. Keyfitz and Kranzer \cite{Ke_Kr_1} worked with a strictly hyperbolic, genuinely nonlinear system derived from a 1-dimensional model for isothermal gas dynamics and observed that there is a large region where the Riemann problem cannot be solved using shocks and rarefactions. They produced approximate unbounded solutions that do not satisfy the equation in the classical weak-solution sense. They also showed that only the first component of the Rankine–Hugoniot relation is satisfied, giving a unique speed $\sigma$ for which any given two states can be joined. Later on, Schecter \cite{Sc} used ideas and methods associated with dynamical systems with geometric flavor (blowing-up approach to geometric singular perturbation problems that lack normal hyperbolicity; see Fenichel \cite{Fe} and Jones \cite{Jo}) to prove the existence of a self-similar viscous solution. See also \cite{Hsu, Ka_Mi, Ka_Mi, Ke, Ke_2, Ke_3, Ke_4, Le_Sl, Ma_Be} and references therein for other solutions involving singular solutions.

\vspace{5mm}

The investigation of singular solutions was mainly focused on when only one state variable
develops the Dirac delta function, and the others are functions with a bounded variation. We have other physically important systems with delta functions in more than one state variable. For example, Mazzotti et al. \cite{Ma_1, Ma_2, Ma_3} numerically studied a model with important applications in modern industry, which exhibits singular solutions arising in two-component chromatography, and both components of the Rankine–Hugoniot relation are not satisfied. Tsikkou \cite{Ts} considered this chromatography system, which exhibits a change of type (hyperbolic and elliptic), performed linear changes in the conserved quantities to obtain a simpler system, and gave a coherent explanation and description of the unbounded solutions.

\vspace{5mm}

It is natural to then ask whether it is possible to predict singular solutions to a system,
find a physical interpretation of their significance, explain the sense in which they satisfy the equation, find a better definition that will describe some broader collection of examples, and check for connections between singular solutions, genuinely nonlinear systems, and change of type (conservation laws which are not everywhere hyperbolic). The model under consideration serves this purpose in addition to the aforementioned physical reasons. From a mathematical point of view, we aim to gain a broader perspective for solving Riemann and Cauchy problems with large data globally using singular solutions as additional building blocks (in possibly generalized schemes). 

\vspace{5mm}
The system (\ref{EQMAIN}) is a special case of 
\begin{equation}\label{EQMODEL}
    \begin{split}\begin{cases}
    \rho_t+\bigg(\rho \Phi\big(\rho, u\big)\bigg)_x&=F\big(\rho, u\big),\\
    \big(\rho u\big)_t+\bigg(\rho u\Phi\big(\rho, u\big)\bigg)_x&=G\big(\rho, u\big),
    \end{cases}
    \end{split}
\end{equation}
where $\Phi(\rho,u)=f(u)-p(\rho)$ is a nonlinear function, and has various applications depending on $\Phi,$ $F,$ and, $G.$ For example, the pressureless Euler system and the macroscopic model for traffic flow by Aw and Rascle \cite{Aw} correspond to $F=G=0, \ \Phi(\rho, u)=u,$ and $F=G=0, \ \ f(u)=u,$ respectively. The literature, see Zhang \cite{Zh_1,Zh_2}, shows that the Riemann problem with pressure laws depending only on the density and $F=0$ has been well studied. Motivated by Li \cite{Li} and references therein, we consider system (\ref{EQMAIN}) with initial data 
\begin{equation}\label{LR_RESTRIC}
    \big(\rho,u\big)\big(x,0\big)=\begin{cases}
    \big(\rho_L,u_L\big)\text{ if } x<0\\\big
    (\rho_R,u_R\big)\text{ if } x>0
\end{cases}
\end{equation}
The paper is organized as follows. In Section II, we use the following substitution:
\begin{align}
    \rho&=v e^{kt}, \ & \ u&=w+\beta t, \ \ & \text{for} \ & \ \eta=k, \label{1.4}\\
     \rho&=v e^{kt}, \ & \ u&=\bigg(w+\frac{\beta}{\eta-k}\bigg)e^{(\eta-k)t}-\frac{\beta}{\eta-k}, \ \ & \text{for} \ & \ \eta\neq k, \label{1.5}
     \end{align}
to transform (\ref{EQMAIN}) to a system of conservation laws and present the numerical method used to verify our analytical results. Section III gives a formal description of the classical Riemann solutions to the system of conservation laws. We use the Rankine-Hugoniot relations to derive the shock curves through a left state and the method of characteristics to get information about the rarefaction curves. All the curves depend on time; therefore, various regions where classical Riemann solutions (using 1-shock, 2-contact discontinuity, and 2-rarefaction) exist evolve in time. On the other hand, the Riemann solution with the right state in Region V consists of delta-shocks. In Section IV, we prove that the singular solution satisfies (\ref{EQMAIN}) in the sense of distributions and we discuss the region time evolution. In Section V, we construct the singular solution to the Riemann problem for the original system (\ref{EQMAIN}), and finally, in Section VI, we present the conclusion.

\section{Preliminaries} 

\subsection {Analysis Preliminaries}

The problem is best split into two cases: $\eta \neq k$ and $\eta=k$. For the $\eta \neq k$ case, the change of variables in (\ref{1.5}) with the restriction $v>0$ is used. With this, we rewrite (\ref{EQMAIN}) into the resulting conservative system

\begin{equation}\label{n_neqk_EQ}
    \begin{split}
    \begin{cases}
        v_t+\bigg[ve^{(\eta-k) t}\bigg(w+\frac{\beta}{\eta-k}\bigg)-v\frac{\beta}{\eta-k}-A\big(ve^{kt}\big)^{\gamma+1}e^{(\eta-k) t}\bigg]_x=0 \\
        \bigg[v\bigg(w+\frac{\beta}{\eta-k}\bigg)\bigg]_t+  \bigg[v \bigg(w+\frac{\beta}{\eta-k}\bigg)^2e^{(\eta-k)t}-  \bigg(w+\frac{\beta}{\eta-k}\bigg)e^{(\eta-k)t}A\big(v e^{kt}\big)^{\gamma+1}-\frac{\beta}{\eta-k}\bigg(w+\frac{\beta}{\eta-k}\bigg)v\bigg]_x=0
    \end{cases}
    \end{split}
\end{equation}
For the other case, $\eta=k$, the transformations in (\ref{1.4}) with restriction $v>0$ are used.  This results in 
\begin{equation} 
\label{n=k_EQ}
    \begin{split}
    \begin{cases}
         v_t+\bigg[v\big(w+\beta t\big)-A\big(ve^{kt}\big)^{\gamma+1}\bigg]_x=0 \\
         \big(vw\big)_t+\bigg[vw\big(w+\beta t\big)-Aw\big(ve^{\eta t}\big)^{\gamma+1}\bigg]_x=0
    \end{cases}
    \end{split}
\end{equation}
Due to (\ref{LR_RESTRIC}), the initial conditions for both systems are
\begin{equation}\label{n_neqkVW}
    \big(v,w\big)\big(x,0\big)=\big(\rho,u\big)\big(x,0\big)=
    \begin{cases}
       \big(v_L,w_L\big)=\big(\rho_L,u_L\big)\text{ if } x<0\\\big(v_R,w_R\big)=\big(\rho_R,u_R\big)\text{ if } x>0
    \end{cases}
\end{equation}

\subsection{Numerical Preliminaries}

The Local Lax-Friedrichs (LLF) scheme was utilized for its simplicity and non-oscillatory behavior. Following the scheme, the spatial and temporal domains were discretized. Neumann boundary conditions were imposed to preserve the left and right states of the solution at the end points. In particular, the $(n+1)^{th}$ temporal solutions were calculated from the neighboring $n^{th}$ solutions and fluxes following from the equation
\begin{align}
    U_{j}^{n+1} = \frac{1}{2}\big(U_{j-1}^n + U_{j+1}^n\big) + \frac{\text{CFL}}{2 \lambda} \big(F_{j+1}^n - F_{j-1}^n\big)
\end{align}
where CFL represents a numerical stability condition for the LLF scheme given by the inequality 
\begin{align}
    \frac{\Delta t}{\Delta x} \lambda \le \frac{1}{2}
\end{align}
in which $\lambda$ represents the maximum wave speed given by the system's eigenvalues. Note that $U$ denotes the vector of unknown variables, and $F$ represents the flux. By construction $\Delta t = \frac{\text{CFL}}{\lambda}$, thus automatically satisfying the CFL condition with $\Delta x = 1$. Also note, as shown by \cite{Tadmor}, satisfying the CFL conditions guarantees that the LLF scheme converges to the physically correct weak solution satisfying entropy conditions.

\vspace{5mm}

To prevent numerical instability, the change of variables $y = vw + v \frac{\beta}{\eta - k}$ was utilized. Thus, the vector of conserved quantities is $H = \begin{bmatrix}
    v & y
\end{bmatrix}^T$ and the flux is written directly in terms of $y$. We then converted $y$ back to $w$ with each iteration. Note that for all numerical figures, $U_L = \begin{bmatrix} v_L & w_L \end{bmatrix}^T \text{ and } U_R = \begin{bmatrix} v_R & w_R \end{bmatrix}^T$.  See also \cite{Lev_1, Lev_2} for additional details on the LLF scheme.

\section{Contact Discontinuity, Shock and Rarefaction}
This section is broken into two cases: $\eta \neq k$ and $\eta = k$. Each case begins by finding the Hugoniot locus, using the Rankine-Hugoniot condition, to find the set of points in state space that may be joined to a fixed left state by a shock satisfying the Lax shock admissibility criterion or by a contact discontinuity. In addition, using the method of characteristics, we show that 1-rarefactions cannot exist and derive information about the 2-rarefactions and how two states can be connected by incorporating one. Next, numerical evidence of the latter is presented and analyzed. Finally, the Hugoniot locus and rarefaction curves are plotted together, splitting the (v, w) state space into regions depending on $\gamma$ and $k$.  Note that numerical evidence is only presented for the $\eta \neq k$ case due to the regions being identical in both cases.  

\subsection{$\eta\neq k$ case}

\subsubsection{Hyperbolicity, Linear Degeneracy and Genuine Nonlinearity}

(\ref{n_neqk_EQ}) is rewritten as
\begin{equation}
    \label{n_neqkEigenE}
    H_t+G_x=0
\end{equation}

where
\begin{equation}
    \label{n_neqkEigenM}
    \begin{cases}
    H=
        \begin{bmatrix}
            v\\ vw+v\frac{\beta}{\eta-k}
        \end{bmatrix}
        \\ \\G=
        \begin{bmatrix}
            ve^{(\eta-k) t}\bigg(w+\frac{\beta}{\eta-k}\bigg)-v\frac{\beta}{\eta-k}-A\big(ve^{kt}\big)^{\gamma+1}e^{(\eta-k) t} \\
            v \bigg(w+\frac{\beta}{\eta-k}\bigg)^2e^{(\eta-k)t}-  \bigg(w+\frac{\beta}{\eta-k}\bigg)e^{(\eta-k)t}A\big(v e^{kt}\big)^{\gamma+1}-\frac{\beta}{\eta-k}\bigg(w+\frac{\beta}{\eta-k}\bigg)v
        \end{bmatrix}
    \end{cases}
\end{equation}
To check whether our system is hyperbolic, we need
\begin{equation}
    \label{n_neqkEigenMT}
    \begin{cases}DH=
        \begin{bmatrix}
            1&0\\w+\frac{\beta}{\eta-k}&v
        \end{bmatrix}\\ \\DG=
        \begin{bmatrix}
            e^{(\eta-k) t}\bigg(w+\frac{\beta}{\eta-k}\bigg)-\frac{\beta}{\eta-k}-A\big(ve^{kt}\big)^{\gamma}e^{(\eta-k) t}\big(\gamma+1\big)&ve^{(\eta-k) t}\\\\
            \bigg(w+\frac{\beta}{\eta-k}\bigg)^2e^{(\eta-k)t}-\frac{\beta}{\eta-k}\bigg(w+\frac{\beta}{\eta-k}\bigg)&
           \\-A\big(\gamma+1\big)\big(ve^{kt}\big)^{\gamma}\bigg(w+\frac{\beta}{\eta-k}\bigg)e^{(\eta-k)t} &-e^{(\eta-k)t}A\big(ve^{kt}\big)^{\gamma+1}-\frac{\beta v}{\eta-k}
        \end{bmatrix}
    \end{cases}
\end{equation}
where $D$ denotes the differential $[\partial/\partial v,\partial/\partial w]$. Solving $\det(DG-\lambda DH)=0$ to obtain the eigenvalues of the system yields 
\begin{equation}
    \label{n_neqkEigenV}
    \begin{cases}
        \lambda_1=\frac{-\beta}{\eta-k}+\bigg[w+\frac{\beta}{\eta-k}-\frac Av\big(ve^{kt}\big)^{\gamma+1}\big(\gamma+1\big)\bigg]e^{(\eta-k)t}
        \\
        \lambda_2=\frac{-\beta}{\eta-k}+\bigg[w+\frac{\beta}{\eta-k}-\frac Av\big(ve^{kt}\big)^{\gamma+1}\bigg]e^{(\eta-k)t}
    \end{cases}
\end{equation}
The corresponding eigenvectors are:
\begin{equation}
    \label{n_neqkEigenVec}
    \begin{cases}r_1=
        \begin{bmatrix}
            1\\0
        \end{bmatrix}
             \\ \\r_2=
        \begin{bmatrix}
            1\\\frac{A\gamma}{v^2}\big(ve^{kt}\big)^{\gamma+1}
        \end{bmatrix}
    \end{cases}
\end{equation}
Here we note that $A<0$ and $\gamma<0$ gives $\lambda_1<\lambda_2$. Furthermore, observe that
\begin{equation}
    \label{n_neqkLinDegen}
    \begin{cases}
        D\lambda_1\cdot r_1=\frac{A}{v^2}\big(\gamma+1\big)\big(ve^{kt}\big)^\gamma e^{\eta t}\big(-v\gamma\big)\neq0
        \\ \\
        D\lambda_2\cdot r_2=\frac{A}{v^2}\big(ve^{kt}\big)^\gamma e^{\eta t}\big(-v\gamma\big)+e^{(\eta-k)t}\frac{A\gamma}{v^2}\big(ve^{kt}\big)^{\gamma+1}=0
    \end{cases}
\end{equation}
Hence, the 1- and 2-characteristic families are genuinely nonlinear and linearly degenerate, respectively. 

\subsubsection{Hugoniot Locus Through a Left State ($v_{-}, w_{-}).$ The Lax Shock Admissibility Criterion}

Let $\sigma(t) = x'(t)$ be the propagation speed.  Using the Rankine-Hugoniot jump conditions, 
\begin{equation}
    \label{RH_relation}
    \begin{cases}
        -\sigma\big(t\big)\big[v\big]_{\text{jump}} + \bigg[ve^{(\eta-k) t}\bigg(w+\frac{\beta}{\eta-k}\bigg)-v\frac{\beta}{\eta-k}-A\big(ve^{kt}\big)^{\gamma+1}e^{(\eta-k) t}\bigg]_{\text{jump}} = 0 \\
        -\sigma\big(t\big)\bigg[vw+\frac{v\beta}{\eta-k}\bigg]_{\text{jump}} + \bigg[v \bigg(w+\frac{\beta}{\eta-k}\bigg)^2e^{(\eta-k)t}- \bigg(w+\frac{\beta}{\eta-k}\bigg)\bigg(e^{(\eta-k)t}A\big(v e^{kt}\big)^{\gamma+1}+\frac{\beta}{\eta-k}v\bigg)\bigg]_{\text{jump}} = 0,
    \end{cases}
\end{equation}   
where $[\cdot]_{\text{jump}}$ denotes the jump across the shock, we conclude that the states that can be connected to $(v_{-},w_{-})$ by a 1-shock or a 2-contact discontinuity lie on the curves
\begin{equation}
    \label{n_neqkShocks1}
        S_1\big(v_{-},w_{-}\big): w=w_{-}
        \end{equation}
or, 
\begin{equation}
\label{n_neqkShocks2}
        C_2 \big(v_{-},w_{-}\big): w=w_{-}-\dfrac{A}{v_{-}}\big(v_{-}e^{kt}\big)^{\gamma+1}+\dfrac{A}{v}\big(v e^{kt}\big)^{\gamma+1},
\end{equation} respectively. These two curves intersect at $(v_{-},w_{-}).$ By (\ref{RH_relation}) we get 
\begin{equation}
    \label{speeds}
\begin{aligned}
\sigma_1\big(t\big)&=\bigg(w_{-}+\dfrac{\beta}{\eta-k}\bigg)e^{(\eta-k)t}-\dfrac{\beta}{\eta-k}-Ae^{(\eta-k)t}e^{k(\gamma+1)t}\ \dfrac{v^{\gamma+1}-v_{-}^{\gamma+1}}{v-v_{-}}\\
\sigma_2\big(t\big)&=\lambda_2\big(v,w\big)=\lambda_2\big(v_{-},w_{-}\big)
\end{aligned}
\end{equation}
For the 1-shock to satisfy the Lax shock admissibility criterion we require 
\begin{equation}
    \label{Lax}
    \lambda_1(v_{-},w_{-})>\sigma_1>\lambda_1(v,w).
\end{equation}
Let $h_1\big(v\big)=v^{\gamma+1}+\gamma v_{-}^{\gamma+1}-\big(\gamma+1\big)v v_{-}^{\gamma}$ and $h_2\big(v\big)=-\gamma v^{\gamma+1}-v_{-}^{\gamma+1}+\big(\gamma+1\big)v^{\gamma} v_{-},$ then (\ref{Lax}) is equivalent to 
\begin{equation}
\begin{cases}
h_1\big(v\big)<0, \ \ h_2\big(v\big)>0 \  \ \text{when} \ \ v>v_{-}\\
h_1\big(v\big)>0, \ \ h_2\big(v\big)<0\ \  \text{when} \ \ v<v_{-},
\end{cases}
\end{equation}
which hold when $$-1<\gamma<0, \ \ \gamma<-1,$$ respectively. This can be easily checked by studying the first and second derivatives of $h_1$ and $h_2.$ Therefore, the admissible parts of the 1-shock curve consist of points with $v>v_{-}$ when $-1<\gamma<0$ and points with $v<v_{-}$ when $\gamma<-1.$

\subsubsection{1-Rarefaction Curve Through a Left State $(v_{-}, w_{-})$}

Recall equations (\ref{n_neqkEigenV}), and (\ref{n_neqk_EQ}). Differentiating (\ref{n_neqk_EQ}) then simplifying with (\ref{n_neqkEigenV}) yields
\begin{equation}
    \begin{split}
        &w_t + \lambda_2 w_x=0\\
        &v_t + \lambda_1 v_x + vw_x e^{(\eta-k)t}=0 
    \end{split}
\end{equation}

We rewrite these equations in matrix form to get
\begin{equation}
    \label{n_neqkRareContraM}
    \begin{bmatrix}
        v_t\\w_t
    \end{bmatrix}
    +
    \begin{bmatrix}
        \lambda_1&ve^{(\eta-k)t}\\0&\lambda_2  
    \end{bmatrix}
    \begin{bmatrix}
        v_x\\w_x
    \end{bmatrix}
        =0
\end{equation}

Note that the eigenvalues and eigenvectors of the matrix $A= \begin{bmatrix}
        \lambda_1&ve^{(\eta-k)t}\\0&\lambda_2  
    \end{bmatrix}$ are the same as (\ref{n_neqkEigenV}) and (\ref{n_neqkEigenVec}).

\begin{figure}[H]
    \centering
    \begin{tikzpicture}[scale=1]
    \begin{axis} [
      ymax=5,
      ymin=0,
      xmax=5,
      xmin=-5,
      axis x line=middle,
      axis y line=middle,
      domain=0:10,
      xtick={-5, -4, ..., 5},
      ytick={0, 1, ..., 5},
      samples=1001,
      yticklabels=\empty,
      xticklabels=\empty,
      xlabel = \(x\),
      ylabel = {\(t\)}
    ]
    \draw[name path = A] (axis cs:-5,2.5) -- (axis cs:0,0);
    \draw[name path = B] (axis cs:-2.5,5) -- (axis cs:0,0);
    \addplot[domain=-5:0,
    dotted,
    very thick
    ] {-x};
    \addplot [black!30] fill between [of = A and B, soft clip={domain=-5:0}];
    \draw[name path = E] (axis cs:1.9,.5) -- (axis cs:1.9,2);
    \draw[name path = C] (axis cs:4.9,2) -- (axis cs:1.9,2);
    \draw[name path = F] (axis cs:4.9,2) -- (axis cs:4.9,.5);
    \draw[name path = D] (axis cs:4.9,.5) -- (axis cs:1.9,.5);
    \addplot [black!10] fill between [of = C and D, soft clip={domain=1.9:4.9}];
    \draw[dotted,very thick] (axis cs:2.2,1.5) -- (axis cs:2.75,1.5);
    \node at (axis cs:3.75,1.5)(B){\color{black}$\frac{dx}{dt}=\lambda_1$};
    \draw[dashed,very thick] (axis cs:2.2,1) -- (axis cs:2.75,1);
    \node at (axis cs:3.75,1)(B){\color{black}$\frac{dx}{dt}=\lambda_2$};
    \node at (axis cs:-4,2.5)(B){\color{black}$R_1$};
    \foreach \t in {1,2,...,5}
           \addplot[domain=-5:5,
            color=black,
            dashed,
             very thick
            ] {x+\t};    
  \end{axis}
\end{tikzpicture}
    \caption{\textit{An illustration of the characteristics for a hypothetical 1-Rarefaction Wave}}
    \label{fig:rarecont}
\end{figure}
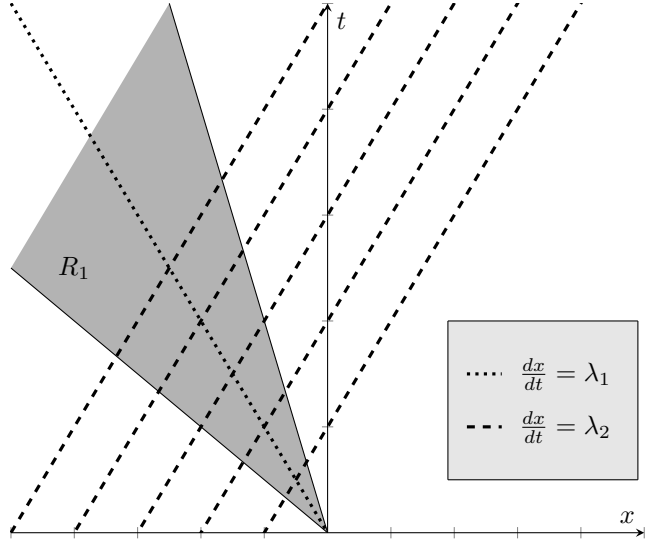

The matrix $A$ is diagonalized in the form $PDP^{-1}$ where 
\begin{equation}\label{D}
    D=
\begin{bmatrix}
    \lambda_1&0\\0&\lambda_2
\end{bmatrix},\end{equation}
\begin{equation}\label{P}
    P=
\begin{bmatrix}
    1&v^2\\0&\gamma A\big(ve^{kt}\big)^{\gamma +1}
\end{bmatrix},\end{equation}
and 
\begin{equation}\label{P-1}
    P^{-1}=
\begin{bmatrix}
     1&\frac{-v^2}{\gamma A\big(ve^{kt}\big)^{\gamma +1}}\\0&\frac1{\gamma A\big(ve^{kt}\big)^{\gamma +1}}
\end{bmatrix}.\end{equation}
These matrices are now used to transform equation (\ref{n_neqkRareContraM}) into
\begin{equation}
    P^{-1}
    \begin{bmatrix}
        v_t\\w_t
    \end{bmatrix}
    +DP^{-1}
    \begin{bmatrix}
        v_x\\w_x
    \end{bmatrix} = 0,
\end{equation}
which yields,
\begin{equation}
\label{n_neqkRareContraF}
    \begin{split}
        &w_t + \lambda_2 w_x =0\\
        &v_t + \lambda_1 v_x -\frac{v^2}{\gamma A \big(v e^{kt}\big)^{\gamma +1}} \big(w_t + \lambda_1 w_x\big)= 0. 
    \end{split}
\end{equation}
This implies that $\frac{dw}{dt}=0$ along 2-characteristics, that is $\frac{dx}{dt}=\lambda_2,$ and thus $w$ is constant. On the other hand, along 1-characteristics where $\frac{dx}{dt}=\lambda_1,$ we get \begin{equation}
    \label{n_neqkRareSlope}
        \frac{dw}{dt} = A \gamma v^{\gamma -1} e^{kt(\gamma+1)} \frac{dv}{dt}.
\end{equation}
As discussed above, when we consider a 1-rarefaction wave connecting the constant left state $(v_-,w_-)$ with another state, see Figure \ref{fig:rarecont} as an illustrative example, $w$ would be constant and equal to $w_-$ in directions given by $\frac{dx}{dt}=\lambda_2.$ 

Therefore, $w$ would stay constant throughout this wave and by (\ref{n_neqkRareSlope}) $v$ would also stay constant since $\frac{dw}{dt}=0$ implies $\frac{dv}{dt}=0.$ By the method of characteristics, we thus conclude that a 1-rarefaction does not exist.


\subsubsection{2-Rarefaction Curve through a left state $(v_{-}, w_{-})$}
We note that for numerical figures the left column displays all 20 iterations, each taking 1000 steps, while the right column displays the latest iteration. Data was renormalized every 100 steps within an error bound of $10^{-7}$ to remove illusory points. Later iterations have a thicker line width. The figures are provided to justify the existence of 2-rarefaction waves numerically. The following sections will discuss the regions and cases mentioned in the figure captions.

\begin{figure}[H]
    \centering
    \includegraphics[width=0.7\linewidth]{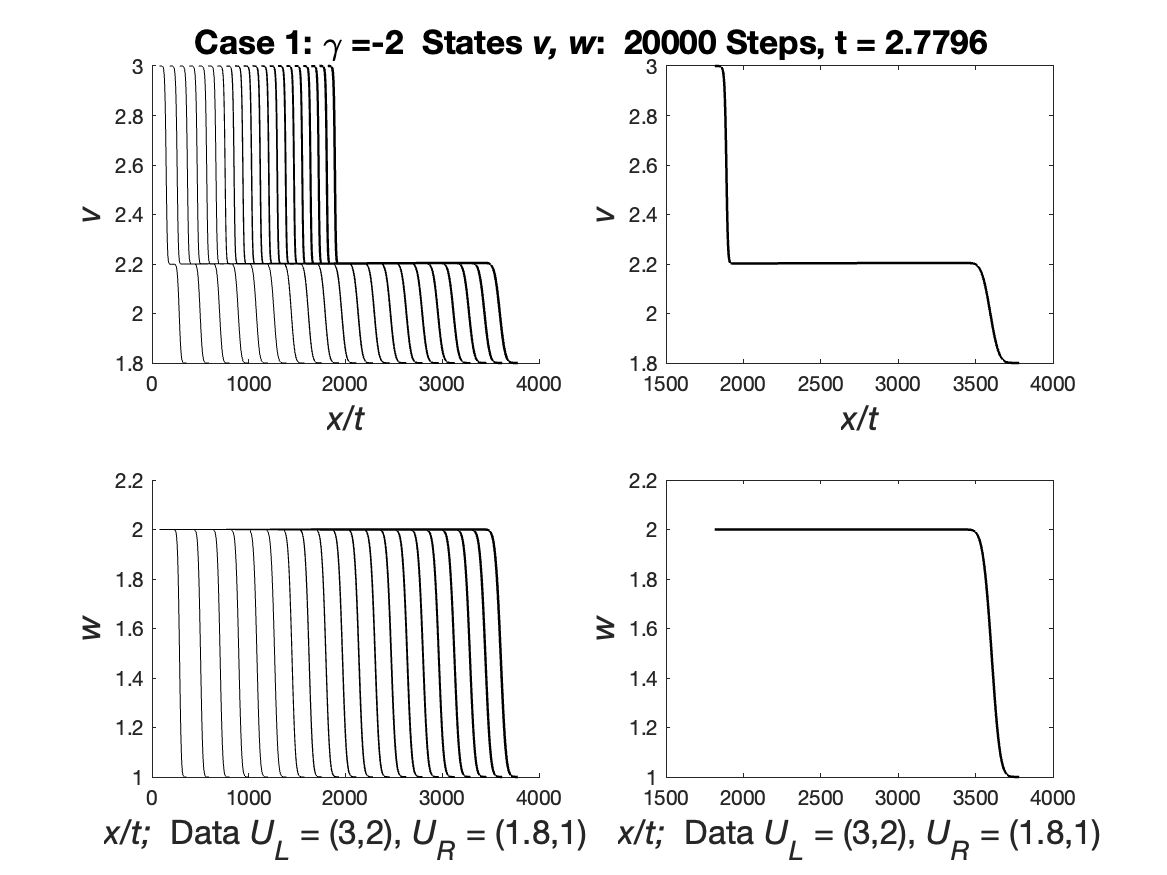}
    \caption{\textit{Region VII in Figure \ref{fig:Regions 4}. \textit{$S_1R_2$}. Parameters: $\gamma = -2, A = -10, \eta = 3, k = 0.01, \beta = 10$.}}
    \label{fig:enter-label}
\end{figure}

\begin{figure}[H]
    \centering
    \includegraphics[width=0.7\linewidth]{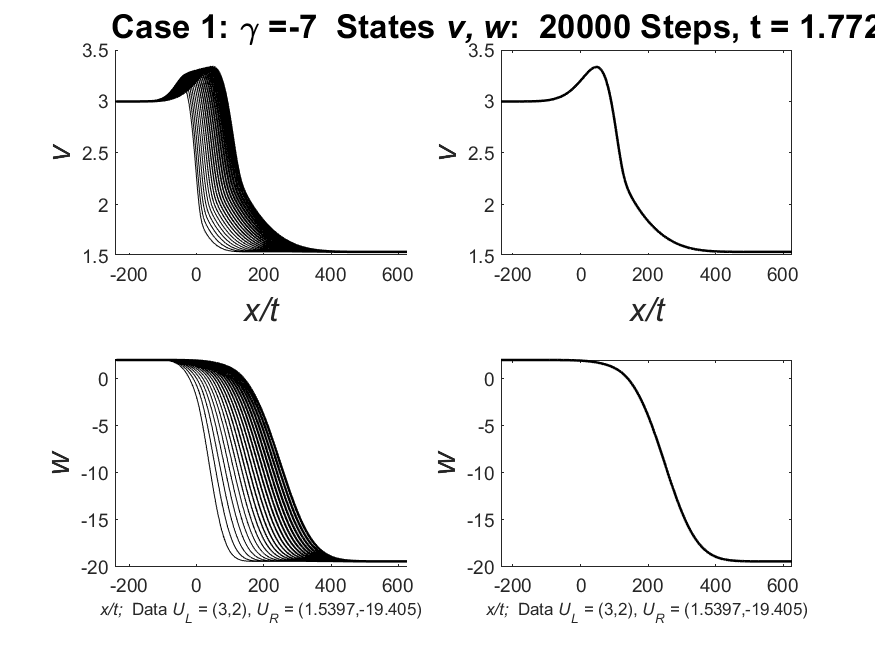}
    \caption{\textit{Region VIII in Figure \ref{fig:Regions 4}. $R_2C_2$. Parameters: $\gamma = -7, A = -500, \eta = 3, k = 0.01, \beta = 10.$}}
    \label{fig:Region VIII Case 1}
\end{figure}

From testing various points, it was found that $R_2$ lies tightly along $C_2$. Numerically, we observed strictly $R_2C_2$ (the Riemann solution consists of a 2-rarefaction followed by a 2-contact discontinuity) when $R_2$ is present for all points tested, but cannot confirm this analytically. Note that the change in $w$ during $R_2$ is extremely small, appearing zero graphically.

\vspace{5mm}

With numerical assurance for a $R_2$ rarefaction, we proceed to analyze its behavior. From (\ref{n_neqkRareSlope}) we can measure the rate of change of $w$ with respect to $t$ on the curve. Since the system is non-autonomous, it is difficult to find the 2-rarefaction curve explicitly. First, we note that a 2-rarefaction $R_2(v_{-},w_{-})$ lies above $C_2(v_{-},w_{-})$; to connect a left state $(v_{-},w_{-})$ with another state $(v,w)$ without crossing characteristics, we require $\lambda_2(v_{-},w_{-})<\lambda_2(v,w)$, ensuring diverging characteristics, as is expected for a rarefaction. To pinpoint the location of $R_2$ we differentiate (\ref{n_neqkShocks2}) along $\frac{dx}{dt}=\lambda_1$ to get
\begin{equation}
\begin{split}
\label{R2}
\frac{dw}{dt}\bigg|_{\text{across} \ C_2}=A\gamma v^{\gamma-1}e^{kt(\gamma+1)}\frac{dv}{dt}+Ak\big(\gamma+1\big)e^{kt(\gamma+1)}\big(v^{\gamma}-v_{-}^{\gamma}\big)
\end{split}
\end{equation}
Note that since a rarefaction wave $R_2$ is a smooth solution we can use (\ref{n_neqkRareSlope}) 
\begin{equation}
    \frac{dw}{dt}\bigg|_{\text{across} \ R_2} =  A \gamma v^{\gamma -1} e^{kt(\gamma+1)} \frac{dv}{dt}
\end{equation}
to study the behavior along 1-characteristics. 
Consequently,
\begin{equation}
\label{dwdt_inequal}
\begin{cases}
\dfrac{dw}{dt}\bigg|_{\text{across} \ C_2} < A\gamma v^{\gamma-1}e^{kt(\gamma+1)}\dfrac{dv}{dt} \ \ \text{when} \ \ k\big(\gamma+1\big)>0  \text{,} \ v<v_{-} \ \ \text{or} \ \ k\big(\gamma+1\big)<0 \text{,} \ v>v_{-}\\ \\ 
\dfrac{dw}{dt}\bigg|_{\text{across} \ C_2} > A\gamma v^{\gamma-1}e^{kt(\gamma+1)}\dfrac{dv}{dt} \ \ \text{when} \ \ k\big(\gamma+1\big)>0 \text{,} \ v>v_{-} \ \ \text{or} \ \ k\big(\gamma+1\big)<0 \text{,} \ v<v_{-}
\end{cases}
\end{equation}
Upon comparing the rates of change of $w$ with respect to $t$ (note that the rates are positive or negative when $v<v_{-}$ or $v>v_{-}$, respectively) on the 2-waves we conclude that when $k(\gamma+1)<0$ the 2-rarefaction curve $R_2$ lies above the 2-contact discontinuity $C_2$ when $v<v_{-}.$ On the other hand, when $k(\gamma+1)>0$ the 2-rarefaction curve $R_2$ lies above $C_2$ when $v>v_{-}.$ 

\vspace{5mm}

Finally, with the assumption that the $R_2$ curve would follow closely above the $C_2$ curve, integration of (\ref{n_neqkRareSlope}) will give a non-explicit equation for $R_2$ as
\begin{equation}
    \label{R2Curveaprrox}
        w= w_{-} +Av^\gamma e^{kt(\gamma +1)} -Av_{-}^\gamma e^{kt_0(\gamma +1)} - \int_{t_0}^t v^\gamma Ak\big(\gamma +1\big) e^{kt(\gamma +1)}dt 
\end{equation}

which is in line with the numerical analysis. These together yield a rough outline of where the $R_2$ curve is located, given by the dotted line in Figures \ref{fig:Regions 4} \ref{fig:Regions 2} \ref{fig:Regions 3} \ref{fig:Regions 1}. This is based off the similarity of (\ref{R2Curveaprrox}) to the equation for $C_2$, as well as the necessary high proximity to $C_2$ to observe $R_2C_2$ instead of $S_1R_2$ numerically. A full derivation of this equation is a topic of future work. 

\vspace{5mm}

All results and analysis of the $R_2$ are expected to hold true for the $\eta = k$ case, due to the inequalities found for the $\eta \neq k $ case reappearing in the former.

\subsubsection{Regions for the Solution of the Riemann Problem}

The curves of our 1-shock $S_1$ and 2-contact discontinuity $C_2$ are given by (\ref{n_neqkShocks1}) and (\ref{n_neqkShocks2}).

The regions are defined further by 
\begin{equation}
    \label{n=kbound}
    S_\delta: w=w_L+\frac{A}{v}\big(ve^{kt}\big)^{\gamma+1},
\end{equation}
representing the limit of the second curve (\ref{n_neqkShocks2}) as $v_L\to\infty$. Additionally, we have 
\begin{equation}
    \label{over}
    S_{o}: w=w_L+\frac{A}{v}\big(ve^{kt}\big)^{\gamma+1}-\frac{A}{v_L}\big(\gamma+1\big)\big(v_Le^{kt}\big)^{\gamma+1}
\end{equation}
representing the max bound of the overcompressive region, which is explained in more detail in section 4.1.1. $R_2$ and all other unknown curves are represented by dotted lines in the $(v, w)$ plane. 

\begin{figure}[H]
\begin{center}
\begin{tikzpicture}
  \begin{axis}[
      ymax=5,
      ymin=-5,
      xmax=5,
      xmin=0,
      axis x line=middle,
      axis y line=left,
      domain=0:10,
      xtick={0, 1, ..., 5},
      ytick={-5, -4, ..., 5},
      samples=1001,
      yticklabels=\empty,
      xticklabels=\empty,
      xlabel = \(v\),
      ylabel = {\(w\)}
    ]
    \addplot [
    domain=.333:5,
    color=black,
    ] {4-(1/x)^2};
    \addplot [
    domain=0:1,
    color=black,
    ] {3};
    \addplot[
    domain=.111:1,
    color=black,
    style=dotted,
    thick,
    ] {5-2/x};
    \addplot [
    domain=.333:5,
    color=black,
    ] {2.75-(1/x)^2};
    \draw[color=black] (axis cs:.4,-1) -- (axis cs:1,-2.5);
    \node at (axis cs:1.25,-2.5)(B){\color{black}VIII};
    \node at (axis cs:2,2)(B){\color{black}V};
    \node at (axis cs:.3,4)(B){\color{black}VI};
    \node at (axis cs:.3,2)(B){\color{black}VII};
    \node at (axis cs:3,4.2)(B){\color{black}$C_2$};
    \node at (axis cs:.5,3.3)(B){\color{black}$S_1$};
    \node at (axis cs:.15,.7)(B){\color{black}$R_2$};
    \node at (axis cs:1.5,2.9)(B){\color{black}$(v_L,w_L)$};
    \node at (axis cs:2,3)(B){\color{black}IX};
    \node at (axis cs:3,3)(B){\color{black}$S_o$};
  \end{axis}
\end{tikzpicture}
\end{center}
\caption{\textit{Regions for $\gamma<-1$ and $k>0$}}
\label{fig:Regions 4}
\end{figure}
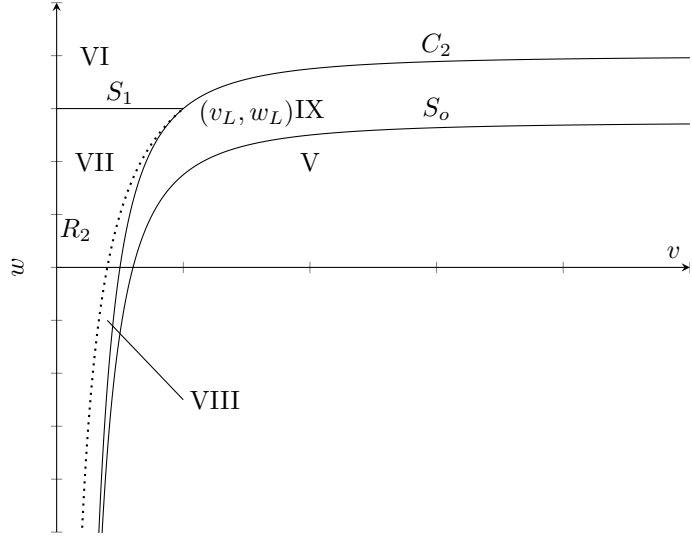

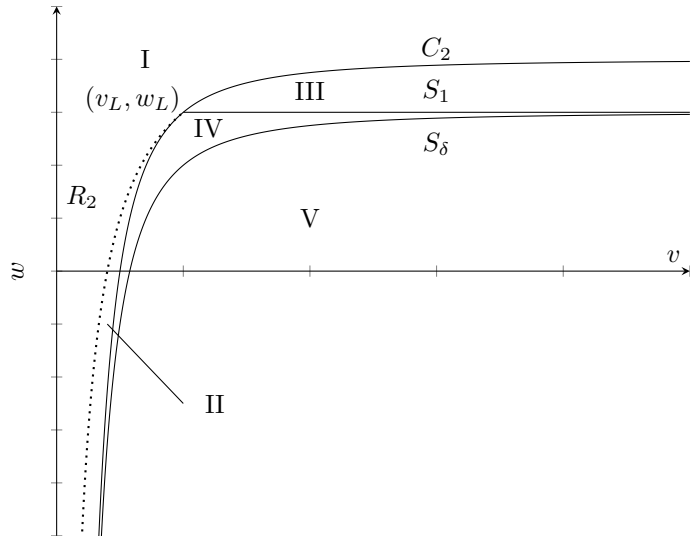
\begin{figure}[H]
\begin{center}
\begin{tikzpicture}
  \begin{axis}[
      ymax=5,
      ymin=-5,
      xmax=5,
      xmin=0,
      axis x line=middle,
      axis y line=left,
      domain=0:10,
      xtick={0, 1, ..., 5},
      ytick={-5, -4, ..., 5},
      samples=1001,
      yticklabels=\empty,
      xticklabels=\empty,
      xlabel = \(v\),
      ylabel = {\(w\)}
    ]
    \addplot [
    domain=.333:5,
    color=black,
    ] {3-(1/x)^2};
    \addplot [
    domain=.333:5,
    color=black,
    ] {4-(1/x)^2};
    \addplot [
    domain=1:5,
    color=black,
    ] {3};
    \addplot[
    domain=.111:1,
    color=black,
    style=dotted,
    thick,
    ] {5-2/x};
    \node at (axis cs:2,3.4)(B){\color{black}III};
    \node at (axis cs:1.2,2.7)(B){\color{black}IV};
    \node at (axis cs:2,1)(B){\color{black}V};
    \node at (axis cs:.7,4)(B){\color{black}I};
    \draw[color=black] (axis cs:.4,-1) -- (axis cs:1,-2.5);
    \node at (axis cs:1.25,-2.5)(B){\color{black}II};
    \node at (axis cs:3,4.2)(B){\color{black}$C_2$};
    \node at (axis cs:3,3.4)(B){\color{black}$S_1$};
    \node at (axis cs:3,2.4)(B){\color{black}$S_\delta$};
    \node at (axis cs:.2,1.4)(B){\color{black}$R_2$};
    \node at (axis cs:.6,3.25)(B){\color{black}$(v_L,w_L)$};
  \end{axis}
\end{tikzpicture}
\end{center}
\caption{\textit{Regions for $-1<\gamma<0$ and $k<0$}}
    \label{fig:Regions 2}
\end{figure}

\begin{figure}[H]
\begin{center}
\begin{tikzpicture}
  \begin{axis}[
      ymax=5,
      ymin=-5,
      xmax=5,
      xmin=0,
      axis x line=middle,
      axis y line=left,
      domain=0:10,
      xtick={0, 1, ..., 5},
      ytick={-5, -4, ..., 5},
      samples=1001,
      yticklabels=\empty,
      xticklabels=\empty,
      xlabel = \(v\),
      ylabel = {\(w\)}
    ]
    \addplot [
    domain=.333:5,
    color=black,
    ] {4-(1/x)^2};
    \addplot [
    domain=0:1,
    color=black,
    ] {3};
    \addplot[
    domain=1:5,
    color=black,
    style=dotted,
    thick,
    ] {-2*exp(1-x)+5};
    \addplot [
    domain=.333:5,
    color=black,
    ] {2.75-(1/x)^2};
    \node at (axis cs:2,2)(B){\color{black}V};
    \node at (axis cs:.3,4)(B){\color{black}VI};
    \node at (axis cs:.3,2)(B){\color{black}VII};
    \node at (axis cs:4,4.5)(B){\color{black}VIII};
    \node at (axis cs:3,4.2)(B){\color{black}$C_2$};
    \node at (axis cs:.5,3.3)(B){\color{black}$S_1$};
    \node at (axis cs:1.3,4.3)(B){\color{black}$R_2$};
    \node at (axis cs:1.5,3)(B){\color{black}$(v_L,w_L)$};
    \node at (axis cs:2,3)(B){\color{black}IX};
    \node at (axis cs:3,3)(B){\color{black}$S_o$};
  \end{axis}
\end{tikzpicture}

\end{center}
\caption{\textit{Regions for $\gamma<-1$ and $k<0$}}
    \label{fig:Regions 3}
\end{figure}
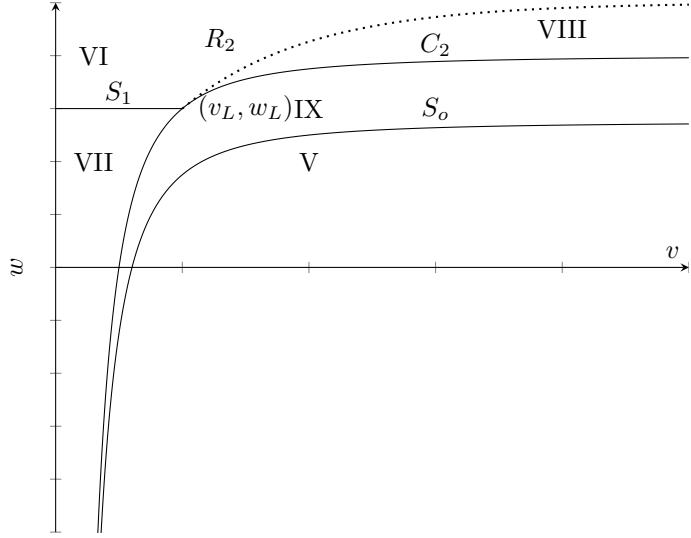

\begin{figure}[H]
\begin{center}
\begin{tikzpicture}
  \begin{axis}[
      ymax=5,
      ymin=-5,
      xmax=5,
      xmin=0,
      axis x line=middle,
      axis y line=left,
      domain=0:10,
      xtick={0, 1, ..., 5},
      ytick={-5, -4, ..., 5},
      samples=1001,
      yticklabels=\empty,
      xticklabels=\empty,
      xlabel = \(v\),
      ylabel = {\(w\)}
    ]
    \addplot [
    domain=.333:5,
    color=black,
    ] {3-(1/x)^2};
    \addplot [
    domain=.333:5,
    color=black,
    ] {4-(1/x)^2};
    \addplot [
    domain=1:5,
    color=black,
    ] {3};
    \addplot[
    domain=1:5,
    color=black,
    style=dotted,
    thick,
    ] {-2*exp(1-x)+5};
    \node at (axis cs:2,3.4)(B){\color{black}III};
    \node at (axis cs:1.2,2.7)(B){\color{black}IV};
    \node at (axis cs:2,1)(B){\color{black}V};
    \node at (axis cs:.5,4)(B){\color{black}I};
    \node at (axis cs: 2.2,4.1)(B){\color{black}II};
    \node at (axis cs:3,4.2)(B){\color{black}$C_2$};
    \node at (axis cs:3,3.4)(B){\color{black}$S_1$};
    \node at (axis cs:3,2.4)(B){\color{black}$S_\delta$};
    \node at (axis cs:1.2,4.1)(B){\color{black}$R_2$};
    \node at (axis cs:.6,3.25)(B){\color{black}$(v_L,w_L)$};
  \end{axis}
\end{tikzpicture}
\end{center}
\caption{\textit{Regions for $-1<\gamma<0$ and $k>0$}}
\label{fig:Regions 1}
\end{figure}
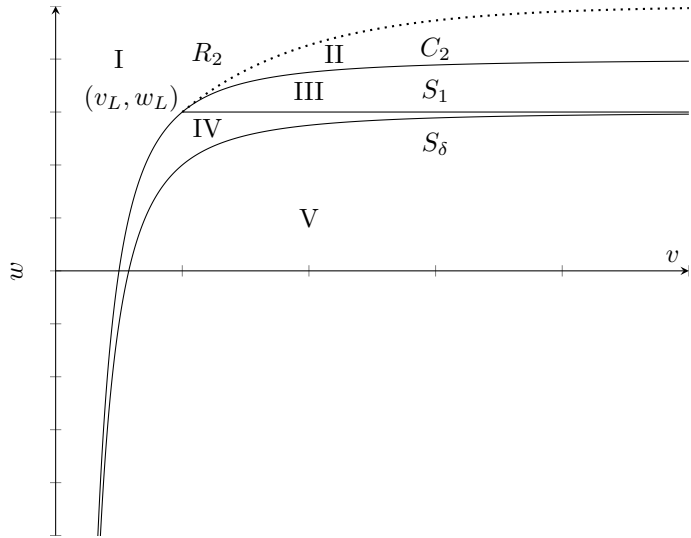

We distinguish four cases: 
\begin{itemize}
\item Case 1, when $\gamma<-1$, $k>0,$ given by Figure \ref{fig:Regions 4},
\item Case 2, when $-1<\gamma<0,k<0$, given by Figure \ref{fig:Regions 2},
\item Case 3, when $\gamma<-1$ and $k<0$ given by Figure \ref{fig:Regions 3},
\item Case 4, when $-1<\gamma<0$ and $k>0$, given by Figure \ref{fig:Regions 1}.
\end{itemize}
For each case, various regions exist that lead to classical and non-classical solutions to the Riemann problem. Specifically, we have
\begin{itemize}
\item a 1-shock followed by a 2-contact discontinuity. The former, given by $x=x_1(t),$ connects $(v_{L},w_{L})$ and a middle state $(v_M,w_M),$ and the latter, given by $x=x_2(t),$ connects the middle state with the right state $(v_R,w_R).$ The middle state can be found explicitly by using (\ref{n_neqkShocks1})-(\ref{n_neqkShocks2}): $$w_M=w_L, \ \ v_M\big(t\big)=\bigg(v_R e^{kt(\gamma+1)}+\frac{w_L-w_R}{A}\bigg)^{1/\gamma} e^{-kt(\gamma+1)/\gamma}.$$ In addition, by (\ref{speeds}) we can find the wave speeds 
\begin{equation}
    \label{speeds_2}
\begin{aligned}
\frac{dx_1}{dt}&=\bigg(w_{L}+\dfrac{\beta}{\eta-k}\bigg)e^{(\eta-k)t}-\dfrac{\beta}{\eta-k}-Ae^{(\eta-k)t}e^{kt(\gamma+1)}\ \dfrac{v_M\big(t\big)^{\gamma+1}-v_L^{\gamma+1}}{v_M\big(t\big)-v_{L}},\\
\frac{dx_2}{dt}&=\frac{-\beta}{\eta-k}+\bigg[w_R+\frac{\beta}{\eta-k}-\frac {A}{v_R}\big(v_Re^{kt}\big)^{\gamma+1}\bigg]e^{(\eta-k)t}.
\end{aligned}
\end{equation}
This solution is possible when the right state is in Region VI (Case 1) or Region VII (Case 3), and Regions III or IV (Cases 2 and 4). The regions are bounded by $S_1(v_L,w_L)$ and $C_2(v_L,w_L)$ or 
$S_{\delta}$ and $C_2(v_L,w_L)$, respectively. 
\item A 1-shock followed by a 2-rarefaction. The solution is possible when the right state is in Region VII (Case 1) or Region VI (Case 3). The regions are bounded by $S_1(v_L,w_L)$, and $R_2(v_L,w_L)$. 
\item A 2-rarefaction followed by a 2-contact discontinuity, which is possible when the right state is in Region I (Cases 2 and 4). The region is bounded by $C_2(v_L,w_L)$, and $R_2(v_L,w_L)$.
\item A 2-rarefaction followed by a 2-contact discontinuity or vice versa. The region is bounded by $C_2(v_L,w_L)$, and $R_2(v_L,w_L)$. This is possible when the right state is in Region VIII (Case 1 and Case 3) or Region II (Case 2 and Case 4). Numerically, we find that the 2-rarefaction comes first in Region VIII in Figure \ref{fig:Regions 4} and expect the same for Region VIII in Figure \ref{fig:Regions 3}. However, due to the growth of $C_2$ into a vertical line during the region shift discussed later, picking a point adequately close to $C_2$ to observe behavior causes the point to almost immediately leave the region, making it impossible to verify. Further work needs to be done to tell if the 2-rarefaction or the 2-contact discontinuity occurs first.
\item A delta-shock which is overcompressible and possible when the right state is in the overcompressive subset of Region V (Cases 1-4). The region is bounded by either $S_{\delta}$ (Cases 2 and 4) or $S_o$ (Cases 1 and 3). Overcompressibility will be discussed in section 4.1.1.
\item Either a delta-shock followed by a 2-wave or a 2-contact discontinuity followed by a delta-shock. This is possible in Region IX (Cases 1 and 3), which is bounded by $C_2(v_L,w_L)$ and $S_{o}$, and in the non-overcompressive subset of Region V in Cases 1-4. More detailed analysis of regions where we expect a combination of a delta-shock and classical solutions will be the subject of future work.

\end{itemize}

\vspace{-0.8cm}

\subsubsection{Numerical Evidence on the Various Regions}
Certain additional restrictions were placed on the numerical constants for the LLF method. We discuss only $A, \eta$ and $k$ here, as $\beta$ has minimal effect on qualitative behavior.

\vspace{5mm}

The magnitude of A is directly related to the size of the regions. We chose $|A| \geq 10$. Note that for large values of $A$, some regions become difficult to access numerically. Due to the $e^{kt}$ term present in $C_2$ and $R_2$, the parameter $k$ controls the rate at which the curves change in time and therefore the regions. In this work, we chose $|k| = 0.01$ or $|k| \geq 0.6$. The former leads to curves that undergo minimal change in time, allowing insight into the initial combinations to enter regions. The latter approximates long-term behavior in time, modeling regional shift over time. While $\eta$ alone does not have a strong effect on region behavior, we require $\eta - k > 0$ due to the factor of $e^{(\eta-k)t}$ in $\lambda_1$ and $\lambda_2$. If $\eta - k < 0$, $\lim_{t\to\infty}e^{(\eta-k)t} = 0$, causing $\lambda_1 = \lambda_2$ in infinite time. We would then lose strict hyperbolicity of our regions, which is incompatible with our implementation and leads to unexpected results. Here, we present the numerical evidence for how the Riemann problem is solved for four cases involving a given left state and a right state in various regions. For Case 1, Region VIII, refer back to Figure \ref{fig:Region VIII Case 1}. We do not show Region VIII (Cases 1 and 3) and Region II (Cases 2 and 4) here, which is reached by $R_2C_2$, since picking an appropriate point is difficult due to the difference between $R_2$ and $C_2$ being extremely small. We identify $S_1$ by it's steep slope and lack of movement in $w$. $R_2$ is classified by a more gradual slope in combination with a ``fanning'' effect over all iterations. $C_2$ has some characteristics of both, often displaying a more gradual slope but with a consistent lack of fanning.  
\vspace{-0.5cm}
\begin{figure}[H]
    \centering
    \includegraphics[width=0.7\linewidth]{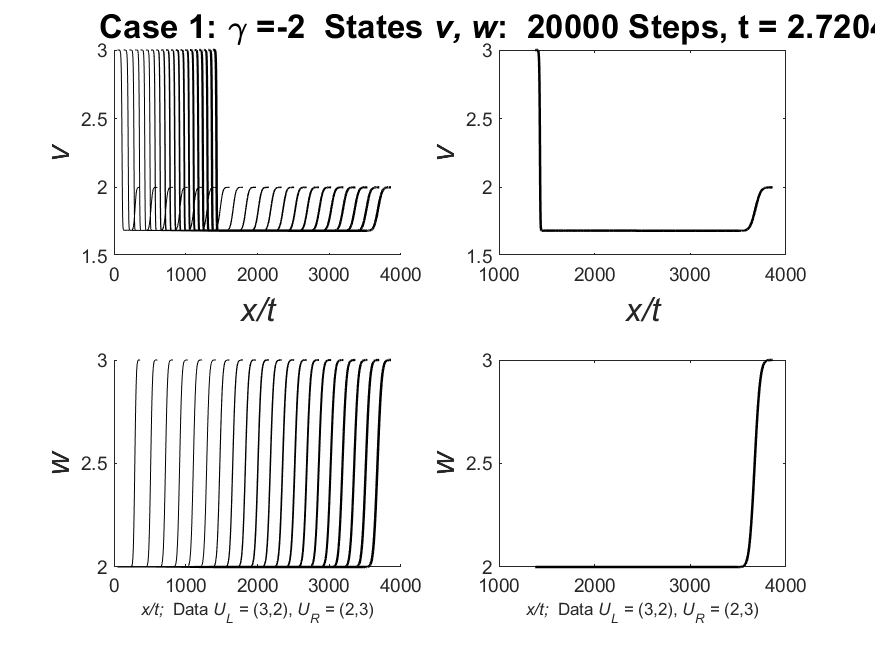}
    \caption{\textit{Region VI of Case 1, $S_1C_2$. Parameters: $\gamma = -2, A = -10, \eta = 3, k = 0.01, \beta = 10$.}}
    \label{fig:Case1_VI_SC}
\end{figure}
\vspace{-0.5cm}
\begin{figure}[H]
    \centering
    \includegraphics[width=0.7\linewidth]{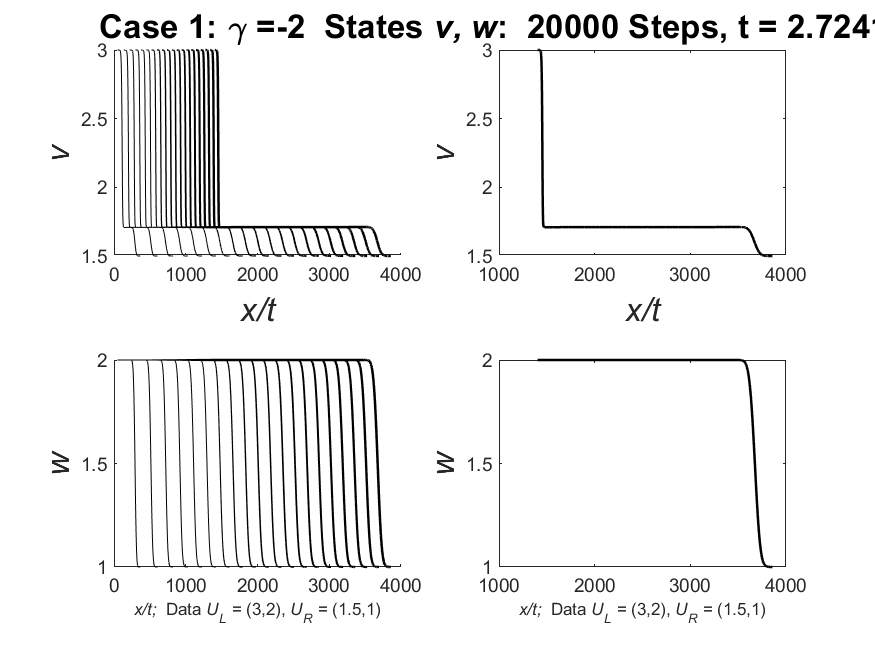}
    \caption{\textit{Region VI of Case 1, $S_1R_2$. Parameters: $\gamma = -2, A = -10, \eta = 3, k = 0.01, \beta = 10$.}}
    \label{fig:Case1_VII_SR}
\end{figure}

\begin{figure}[H]
    \centering
    \includegraphics[width=0.7\linewidth]{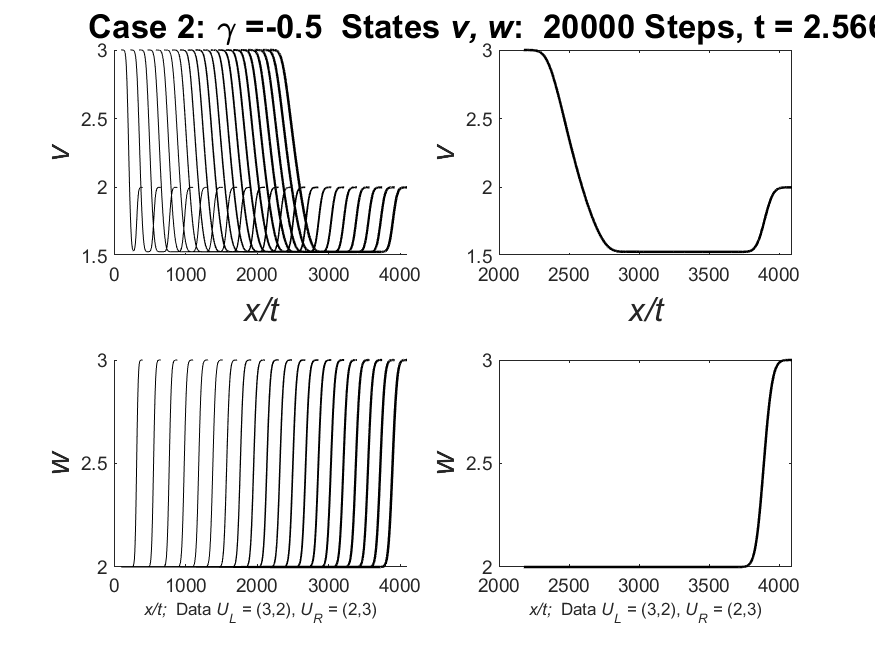}
    \caption{\textit{Region I of Case 2, $R_2C_2$. Parameters: $\gamma = -0.5, A = -10, \eta = 3, k = -0.01, \beta = 10$.}}
    \label{fig:Case2_I_RC}
\end{figure}

\begin{figure}[H]
    \centering
    \includegraphics[width=0.7\linewidth]{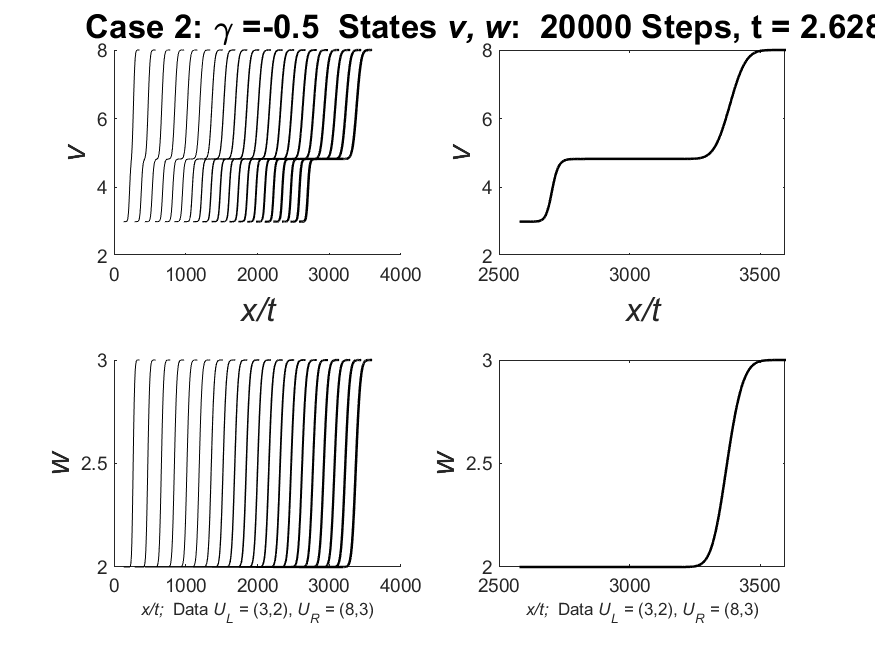}
    \caption{\textit{Region I of Case 2, $S_1C_2$. Parameters: $\gamma = -0.5, A = -10, \eta = 3, k = -0.01, \beta = 10$.}}
    \label{fig:Case2_III_SC}
\end{figure}

\begin{figure}[H]
    \centering
    \includegraphics[width=0.7\linewidth]{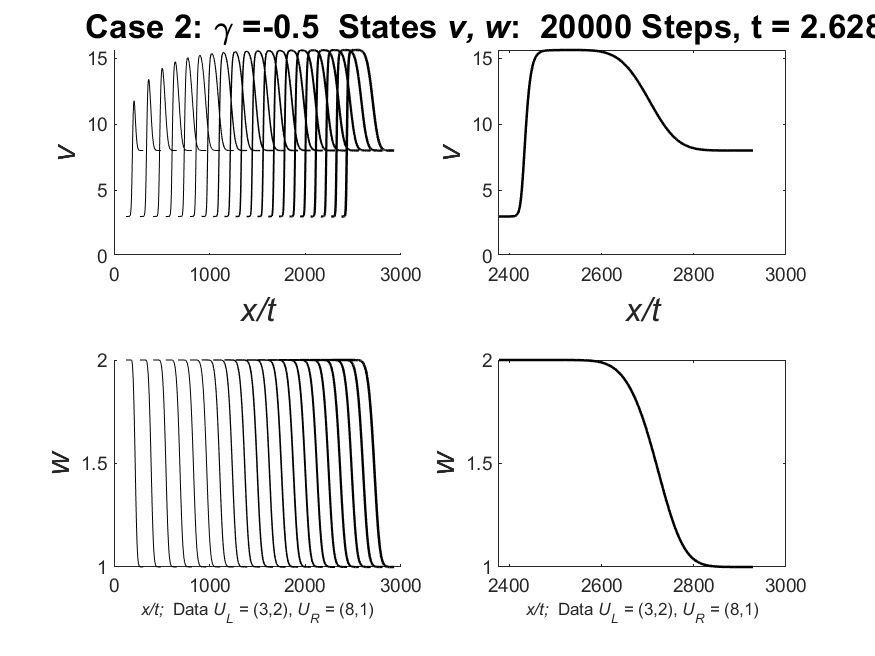}
    \caption{\textit{Region I of Case 2, $S_1C_2$. Parameters: $\gamma = -0.5, A = -10, \eta = 3, k = -0.01, \beta = 10$.}}
    \label{fig:Case2_IV_SC}
\end{figure}

\begin{figure}[H]
    \centering
    \includegraphics[width=0.7\linewidth]{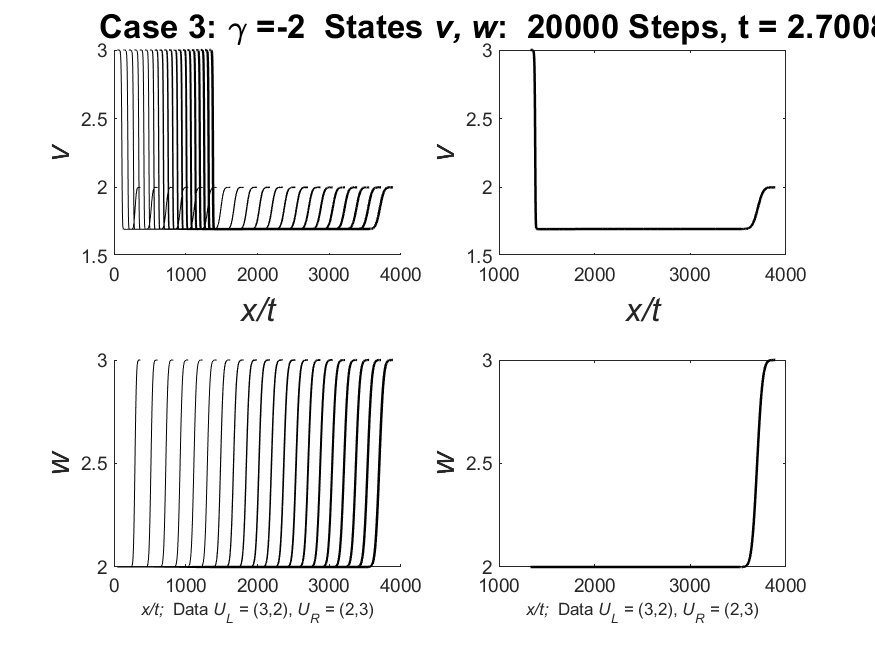}
    \caption{\textit{Region VI of Case 3, $S_1R_2$. Parameters: $\gamma = -2, A = -10, \eta = 3, k = -0.01, \beta = 10$.}}
    \label{fig:Case3_VI_SR}
\end{figure}

\begin{figure}[H]
    \centering
    \includegraphics[width=0.7\linewidth]{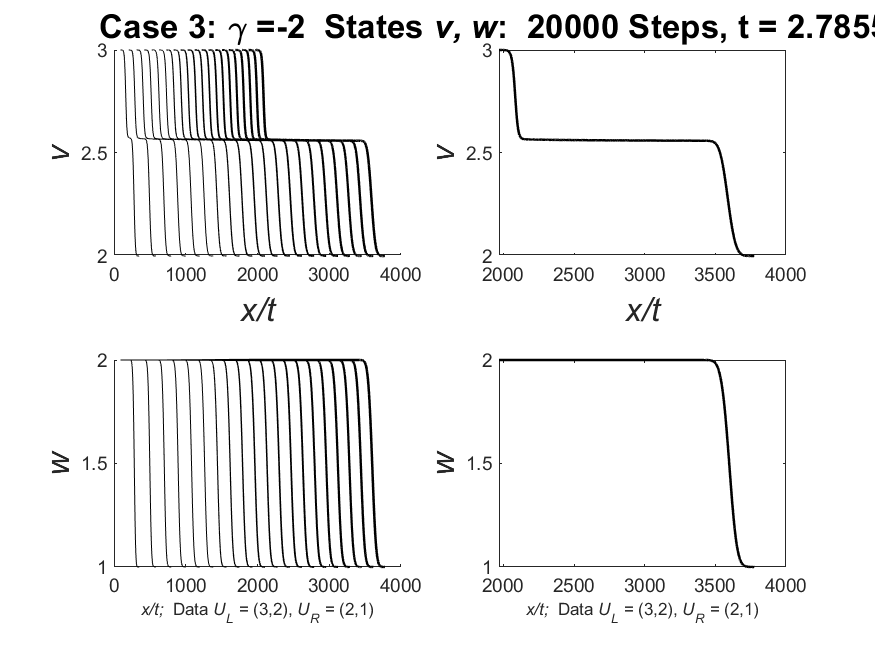}
    \caption{\textit{Region VII of Case 3, $S_1C_2$. Parameters: $\gamma = -2, A = -10, \eta = 3, k = -0.01, \beta = 10$.}}
    \label{fig:Case3_VII_SC}
\end{figure}

\begin{figure}[H]
    \centering
    \includegraphics[width=0.7\linewidth]{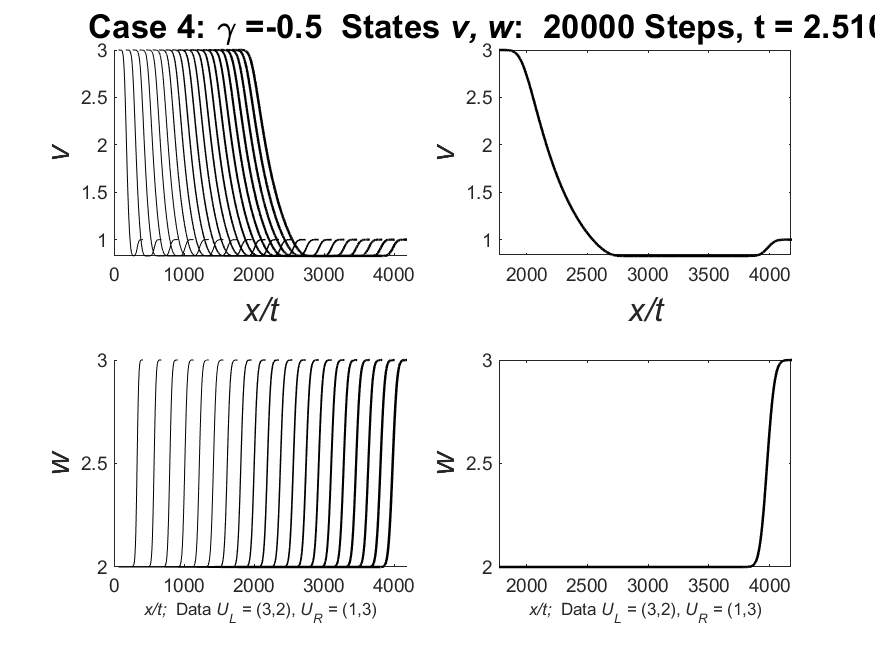}
    \caption{\textit{Region I of Case 4, $R_2C_2$. Parameters: $\gamma = -0.5, A = -10, \eta = 3, k = 0.01, \beta = 10$.}}
    \label{fig:Case4_I_RS}
\end{figure}

\begin{figure}[H]
    \centering
    \includegraphics[width=0.7\linewidth]{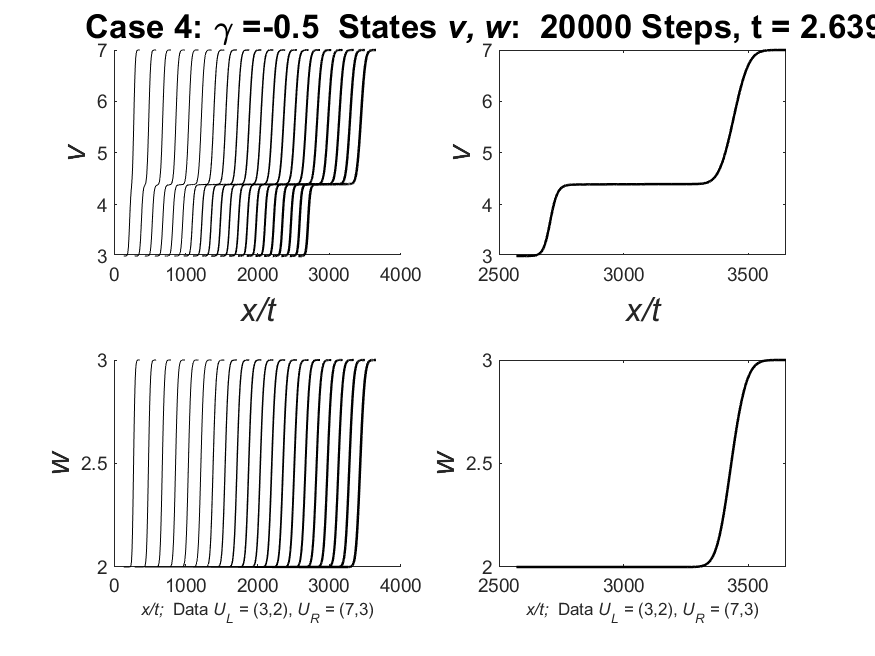}
    \caption{\textit{Region III of Case 4, $S_1C_2$. Parameters: $\gamma = -0.5, A = -10, \eta = 3, k = 0.01, \beta = 10$.}}
    \label{fig:Case4_III_SC}
\end{figure}

\begin{figure}[H]
    \centering
    \includegraphics[width=0.7\linewidth]{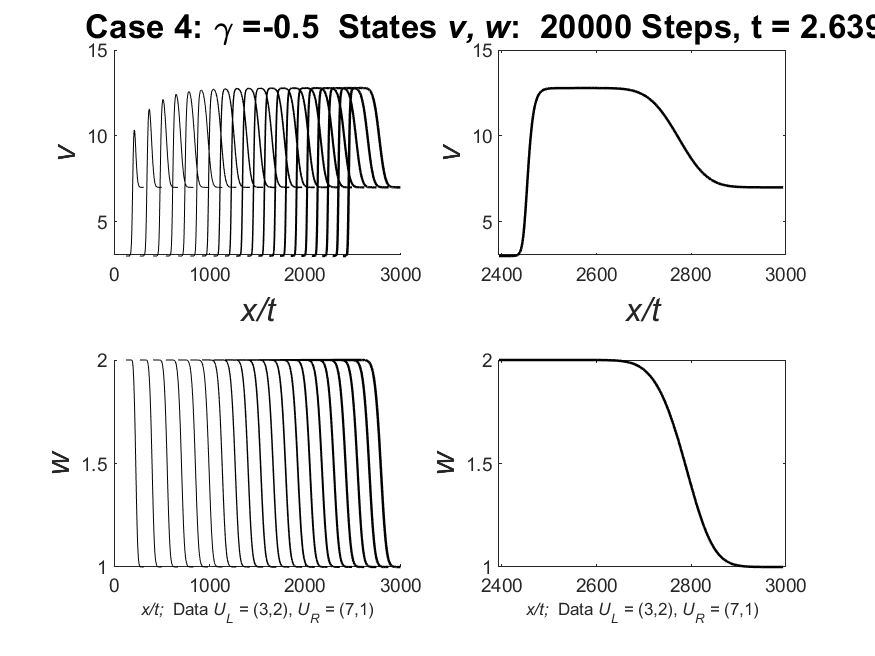}
    \caption{\textit{Region III of Case 4, $S_1C_2$. Parameters: $\gamma = -0.5, A = -10, \eta = 3, k = 0.01, \beta = 10$.}}
    \label{fig:Case4_IV_SC}
\end{figure}

\subsection{$\eta = k$ Case}

\subsubsection{Hyperbolicity, Linear Degeneracy and Genuine Nonlinearity}

Equation (\ref{n=k_EQ}) can be rewritten as 
\begin{equation}
    \label{n=k_EQ_MAT}
    H_t + G_x = 0 
\end{equation}

where $G$ and $H$ are taken to be
\begin{equation}
    \label{n=kG+H}
    \begin{split}
        \begin{cases}
        H=
            \begin{bmatrix}
                v\\vw
            \end{bmatrix}\\ \\
            G=
            \begin{bmatrix}
                v\big(w+\beta t\big)-A\big(ve^{kt}\big)^{\gamma+1}\\ vw\big(w+\beta t\big)- Aw\big(ve^{kt}\big)^{\gamma+1}
            \end{bmatrix}
        \end{cases}
    \end{split}
\end{equation}

To check, once again, whether our system is hyperbolic, $\det(DG-\lambda DH)=0$ is solved, to find the eigenvalues:
\begin{equation}
    \label{n=k_EVAL}
    \begin{split}
        \begin{cases}
            \lambda_1=\big(w+\beta t\big)-A\big(\gamma+1\big)v^\gamma e^{kt(\gamma+1)} \\
            \lambda_2=\big(w+\beta t\big)-Av^\gamma e^{kt(\gamma+1)}
        \end{cases}
    \end{split}
\end{equation}

with the corresponding eigenvectors:
\begin{equation}
    \label{n=k_EVEC}
    \begin{split}
        \begin{cases}
            r_1=
            \begin{bmatrix}
                1\\0
            \end{bmatrix}
            \\ \\
            r_2=
            \begin{bmatrix}
                1\\ A\gamma v^{\gamma-1} e^{kt(\gamma+1)}
            \end{bmatrix}
        \end{cases} 
    \end{split}
\end{equation}

Additionally, 
\begin{equation}
    \label{n=kLINDEG}
    \begin{split}
        \begin{cases}
            D\lambda_1\cdot r_1=-A\gamma\big(\gamma+1\big)v^{\gamma-1} e^{kt(\gamma+1)}\neq 0 \\
            D\lambda_2\cdot r_2=-A\gamma v^{\gamma} e^{kt(\gamma+1)}+A\gamma v^{\gamma-1} e^{kt(\gamma+1)}=0
        \end{cases}
    \end{split}
\end{equation}

Similar to the $\eta \neq k$ case, the 1- and 2-characteristic families are genuinely nonlinear and linearly degenerate, respectively. 

\subsubsection{Hugoniot Locus Through a Left State $(v_{-},w_{-})$. The Lax Shock Admissibility Criterion}

Once again, the Rankine-Hugoniot jump conditions (\ref{RH_relation}) are checked, 
\begin{equation}
    \label{n_eqk_shk_ot}
    \begin{cases}
        -\sigma\big(t\big)[v]_{\text{jump}}+\bigg[v\big(w+\beta t\big)-A\big(v e^{kt}\big)^{\gamma+1}\bigg]_{\text{jump}}&=0\\
        -\sigma\big(t\big)[vw]_{\text{jump}}+\bigg[vw\big(w+\beta t\big) -Aw\big(ve^{kt}\big)^{\gamma+1} \bigg]_{\text{jump}}&=0
    \end{cases}
\end{equation}
resulting in
\begin{equation}
    \label{n_eqk_shk_oo}
    \begin{cases}
        S_1\big(v_{-},w_{-}\big): w=w_{-}\\ \\
        C_2\big(v_{-},w_{-}\big): w=w_{-}-\dfrac{A}{v_{-}}\big(v_{-}e^{kt}\big)^{\gamma+1}+\dfrac{A}{v}\big(v e^{kt}\big)^{\gamma+1}.
    \end{cases}
\end{equation}
Therefore, the states that can be connected to $(v_{-},w_{-})$ by a 1-shock or a 2-contact discontinuity lie on the curves (\ref{n_eqk_shk_oo}). By (\ref{n_eqk_shk_ot})

\begin{align}\label{speedsn=k}
    \sigma_1\big(t\big)&=w_{-}+\beta t-A e^{kt(\gamma+1)}\           \dfrac{v^{\gamma+1}-v_{-}^{\gamma+1}}{v-v_{-}}\\
    \sigma_2\big(t\big)&=\lambda_2\big(v,w\big)=\lambda_2\big(v_{-},w_{-}\big)
\end{align}

are obtained. Again, (\ref{Lax}) is checked to ensure that the Lax shock admissibility criterion is satisfied. The first and second inequalities give
\begin{equation}
    \label{n_neqk_EC_EQT}
    v_{-}^{\gamma}\big(\gamma+1\big)>\frac{v^{\gamma+1}-v_{-}^{\gamma+1}}{v-v_{-}}
\end{equation}
and
\begin{equation}
    \label{n_neqk_EC_EQO}
    v^{\gamma}\big(\gamma+1\big)<\frac{v^{\gamma+1}-v_-^{\gamma+1}}{v-v_-},
\end{equation}
respectively. Both are equivalent to the corresponding inequalities for the case $\eta\neq k$. Thus, $S_1$ exists for $v>v_{-}$ if $-1<\gamma<0$ and for $v<v_{-}$ if $\gamma<-1$. The rest of the analysis (non-existence of 1-rarefactions, existence and location of 2-rarefactions, regions, etc.) is identical to the $n\neq k$ case and will be omitted for brevity.

\section{Delta-Shocks}
In the overcompressibe subset of Region V, there is no solution that is piecewise smooth, and bounded. Therefore, in order to establish existence in a space of measures from a mathematical perspective, a solution containing a weighted $\delta$-measure (or $\delta$-shock) supported on a curve needs to be constructed. For these singular solutions, we need to consider physical constraints. If $w$ becomes unbounded, the system's velocity must approach infinity. However, this scenario is not physically possible, as we know that the speed of light bounds the velocity of all particles in the universe. Furthermore, it is also not physically feasible for both $v$ and $w$ to be unbounded, as it implies a finite amount of mass suddenly becoming infinite and then returning to a finite value. Therefore, the only possible case is for $v$ to be unbounded. Assuming infinite density in situations where the fluid volume is nearly infinitesimal is physically reasonable. We have also observed numerically the presence of the Dirac delta measure in $v$ only.
\subsection{$\eta \neq k$ Case}
We define a two-dimensional weighted $\delta$-measure $\omega(s)\delta_S$ supported on a smooth curve $S=\{(x(s),t(s)): c \leq s \leq d\}$ by 
$$\bigg\langle \omega\big(\cdot\big)\delta_S,\psi\big(\cdot,\cdot\big)\bigg\rangle=\int_a^b \omega\big(t\big(s\big)\big)\psi\big(x\big(s\big),t\big(s\big)\big)\ ds$$ for all $\psi\in C_0^{\infty} \big(\mathbb{R}\times\mathbb{R}^+\big).$

Following the above reasoning, the definition of solutions in the sense of distributions is as follows.

\begin{figure}[H]
    \centering
    \includegraphics[width=.6\textwidth]{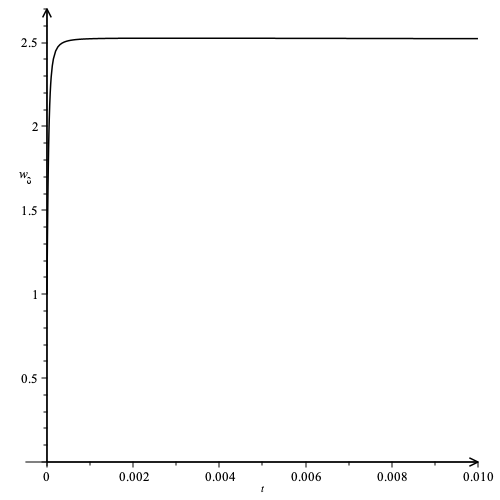}
    \caption{\textit{Numerical solution for $\omega_\delta$ when $\eta=-k\gamma$}}
    \label{fig:Delta Shock}
\end{figure}

{\bf Definition}: A pair $(v, w)$ are known as a \textit{delta-shock type solution} to the system with Riemann data in the sense of distributions if there exists a smooth curve $S=\{(x(t),t: 0 \leq t <\infty\}$ and a weight $\omega_1 \in C^{1}(S) $ such that $v$ and $w$ are represented in the following way
\begin{equation}
    \label{n=kDeltaDef}
    \begin{aligned}
    \big(v,w\big)\big(x,t\big) &= \bigg(v_0\big(x,t\big)+\omega_1\big(t\big)\delta_S,w_0\big(x,t\big)\bigg)\\
    &=
    \begin{cases}
        \big(v_L,w_L\big) \text{,     } & x<x\big(t\big)
        \\
        \bigg(v_{\delta}\big(t\big)+\omega_1\big(t\big)\delta\big(x-x\big(t\big)\big), w_{\delta}\big(t\big)\bigg) \text{,  } & x=x\big(t\big)
        \\
        \big(v_R,w_R\big) \text{,    } & x>x\big(t\big). \text{ }
    \end{cases}
    \end{aligned}
\end{equation}
where $\delta(\cdot)$ is the standard Dirac measure, (therefore $w$ is $v$-measurable, and for example, $v\bigg(w+\frac{\beta}{\eta-k}\bigg)$ can be understood as a Radon measure) and satisfy (\ref{n_neqk_EQ}) in the sense of distributions:
\begin{equation}
        \bigg\langle v,\phi_t\bigg\rangle+\bigg\langle v\bigg(w+\frac{\beta}{\eta-k}\bigg)e^{(\eta-k)t}-\frac{v\beta}{\eta-k}-A\big(ve^{kt}\big)^{\gamma+1}e^{(\eta-k)t},\phi_x\bigg\rangle=0
        \label{Rankine_n_neqk}
\end{equation}
\begin{equation}
        \bigg\langle v\bigg(w+\frac{\beta}{\eta-k}\bigg), \phi_t \bigg\rangle + \bigg\langle v \bigg(w+\frac{\beta}{\eta-k}\bigg)^2e^{(\eta-k)t}-\bigg(w+\frac{\beta}{\eta-k}\bigg)e^{(\eta-k)t}A\big(v e^{kt}\big)^{\gamma+1}-\frac{\beta}{\eta-k}\bigg(w+\frac{\beta}{\eta-k}\bigg)v,\phi_x\bigg\rangle = 0 \label{Rankine_n_neqk_2}
\end{equation}
for every $\phi \in {C_0}^{\infty} \big(\mathbb{R} \times \mathbb{R}^+\big),$ where
\begin{equation}
\notag
\begin{aligned}
\bigg\langle v,\psi\bigg\rangle&=\int_0^{\infty}\int_{-\infty}^{\infty} v_0 \psi \ dx \ dt+\bigg\langle \omega_1(t)\delta_S,\psi\bigg\rangle,\\
\bigg\langle v\bigg(w+\frac{\beta}{\eta-k}\bigg),\psi\bigg\rangle&=\int_0^{\infty}\int_{-\infty}^{\infty} v_0\bigg(w_0+\frac{\beta}{\eta-k}\bigg) \psi \ dx \ dt+\bigg\langle \omega_1(t)\bigg(w_{\delta}(t)+\frac{\beta}{\eta-k}\bigg)\delta_S,\psi\bigg\rangle,\\
\bigg\langle A\big(ve^{kt}\big)^{\gamma+1}e^{(\eta-k)t},\psi\bigg\rangle&=\int_0^{\infty}\int_{-\infty}^{\infty} A\big(v_0 e^{kt}\big)^{\gamma+1}e^{(\eta-k)t} \psi \ dx \ dt,\\
\end{aligned}
\end{equation}
since $\gamma<0.$ The remaining integrals in (\ref{Rankine_n_neqk})-(\ref{Rankine_n_neqk_2}) are similar.Therefore (\ref{Rankine_n_neqk})-(\ref{Rankine_n_neqk_2}) give
\begin{equation}
    \label{Rankine_n_neqk3}
    \begin{split}
        &\int_0^\infty\int_{-\infty}^{x(t)}\bigg(v_L\phi_t+\bigg(v_L\bigg(w_L+\frac{\beta}{\eta-k}\bigg)e^{(\eta-k)t}-v_L\frac{\beta}{\eta-k}-Av_L^{\gamma+1}e^{k(\gamma+1)t}e^{(\eta-k)t}\bigg)\phi_x\bigg)dxdt \\&+\int_0^\infty\int_{x(t)}^\infty \bigg(v_R\phi_t+\bigg(v_R\bigg(w_R+\frac{\beta}{\eta-k}\bigg)e^{(\eta-k)t}-v_R\frac{\beta}{\eta-k}-Av_R^{\gamma+1}e^{k(\gamma+1)t}e^{(\eta-k)t}\bigg)\phi_x\bigg)dxdt
        \\&+\int_0^\infty\bigg( \omega_1\phi_t+\bigg(\omega_1\bigg(w_{\delta}+\frac{\beta}{\eta-k}\bigg)e^{(\eta-k)t}-\omega_1\frac{\beta}{\eta-k}\bigg)\phi_x\bigg)dt\\&=0
    \end{split}
\end{equation}
and 
\begin{equation}
    \label{Rankine_n_neqk4}
    \begin{split}
&\int_0^\infty\int_{-\infty}^{x(t)}\bigg(v_L\bigg(w_L+\frac{\beta}{\eta-k}\bigg)\phi_t+\bigg(v_L\bigg(w_L+\frac{\beta}{\eta-k}\bigg)^2e^{(\eta-k)t}-\frac{\beta}{\eta-k}v_L\bigg(w_L+\frac{\beta}{\eta-k}\bigg) \\ &-Av_L^{\gamma+1}e^{k(\gamma+1)t}\bigg(w_L+\frac{\beta}{\eta-k}\bigg)e^{(\eta-k)t}\bigg)\phi_x\bigg)dxdt\\&+\int_0^\infty\int_{x(t)}^\infty \bigg(v_R\bigg(w_R+\frac{\beta}{\eta-k}\bigg)\phi_t+\bigg(v_R\bigg(w_R+\frac{\beta}{\eta-k}\bigg)^2e^{(\eta-k)t}-\frac{\beta}{\eta-k}v_R\bigg(w_R+\frac{\beta}{\eta-k}\bigg)\\&-Av_R^{\gamma+1}e^{k(\gamma+1)t}\bigg(w_R+\frac{\beta}{\eta-k}\bigg)e^{(\eta-k)t}\bigg)\phi_x\bigg)dxdt\\&+\int_0^\infty \bigg(\omega_1\bigg(w_{\delta}+\frac{\beta}{\eta-k}\bigg)\phi_t+\bigg(\frac{\beta}{\eta-k}\bigg)e^{(\eta-k)t}+\omega_1\bigg(w_{\delta}+\frac{\beta}{\eta-k}\bigg)^2e^{(\eta-k)t}\\&-\frac{\beta}{\eta-k}\bigg(\omega_1\bigg(w_{\delta}+\frac{\beta}{\eta-k}\bigg)\bigg)-Av_{\delta}^{\gamma+1}e^{k(\gamma+1)t}\omega_2e^{(\eta-k)t}\bigg)\phi_x\bigg)dt\\
&=0,
    \end{split}
\end{equation}
respectively. To be able to integrate along $x=x(t)$, we require
\begin{equation}
    \label{shock}
\frac{dx\big(t\big)}{dt} = \sigma\big(t\big) = \bigg(\omega_\delta\big(t\big)+\frac{\beta}{\eta-k}\bigg)e^{(\eta-k)t}-\frac{\beta}{\eta-k}.
\end{equation} 
We apply Green's Theorem, to write (\ref{Rankine_n_neqk3})-(\ref{Rankine_n_neqk4}) as
\begin{equation}\label{Rankine_n_neqk5}
\begin{split}&\int_0^\infty \bigg(-\big[v\big]_{\text{jump}}\sigma+\bigg[v\bigg(w+\frac{\beta}{\eta-k}\bigg)e^{(\eta-k)t}-\frac{v\beta}{\eta-k}-A\big(ve^{kt}\big)^{\gamma+1}e^{(\eta-k)t}\bigg]_{\text{jump}} -\frac{d\omega_1}{dt}\bigg)\phi dt= 0 \end{split}\end{equation}
and 
\begin{equation}\label{Rankine_n_neqk6}\begin{split}&\int_0^\infty \bigg(-\bigg[v\bigg(w+\frac{\beta}{\eta-k}\bigg)\bigg]_{\text{jump}}\sigma+\bigg[\bigg(w+\frac{\beta}{\eta-k}\bigg)^2ve^{(\eta-k)t}-\bigg(w+\frac{\beta}{\eta-k}\bigg)e^{(\eta-k)t}A\big(ve^{kt}\big)^{\gamma+1}\\&-\frac{\beta}{\eta-k}\bigg(w+\frac{\beta}{\eta-k}\bigg)v\bigg]_{\text{jump}}-\frac{d}{dt}\bigg(\omega_1\bigg(\omega_\delta+\frac{\beta}{\eta-k}\bigg)\bigg)\bigg)\phi dt= 0\end{split}\end{equation}
where $[\ \cdot \ ]_{\text{jump}}=\cdot_L-\cdot_R.$ Thus, if we also require 
\begin{align}
            \frac{d\omega_1}{dt}=& -\big[v\big]_{\text{jump}}\sigma+\bigg[v\bigg(w+\frac{\beta}{\eta-k}\bigg)e^{(\eta-k)t}-\frac{v\beta}{\eta-k}-A\big(ve^{kt}\big)^{\gamma+1}e^{(\eta-k)t}\bigg]_{\text{jump}}
         \label{n=kRanHugRel_1}\\ \label{n=kRanHugRel_2}
            \frac{d}{dt}\bigg(\omega_1\bigg(\omega_\delta+\frac{\beta}{\eta-k}\bigg)\bigg)=& -\bigg[v\bigg(w+\frac{\beta}{\eta-k}\bigg)\bigg]_{\text{jump}}\sigma+\bigg[\bigg(w+\frac{\beta}{\eta-k}\bigg)^2ve^{(\eta-k)t}
            \\ \notag &-\bigg(w+\frac{\beta}{\eta-k}\bigg)e^{(\eta-k)t}A\big(ve^{kt}\big)^{\gamma+1}-\frac{\beta}{\eta-k}\bigg(w+\frac{\beta}{\eta-k}\bigg)v\bigg]_{\text{jump}} 
\end{align}
then $(v,w)$ satisfies the system in the sense of distributions, that is (\ref{Rankine_n_neqk})-(\ref{Rankine_n_neqk_2}) hold for every test function $\phi.$ 
It should be noted that there is a Rankine-Hugoniot deficit in both components due to (\ref{n=kRanHugRel_1})-(\ref{n=kRanHugRel_2}), just like in the chromatography model by Mazzotti et al. \cite{Ma_1, Ma_2, Ma_3}. To get more information about $\omega_{1},$ $w_{\delta}$ and $x(t),$ (with $x(0)=0,$ $\omega_1(0)=0$) we return to the original variables $\rho_L(=v_L)$, $\rho_R(=v_R),$ $u_L(=w_L)$, $u_R(=w_R),$ and substitute (\ref{shock}) into (\ref{n=kRanHugRel_1}) to get 
\begin{align}
            \frac{d\omega_1}{dt}=& -w_{\delta}e^{(\eta-k)t}\big[\rho\big]_{\text{jump}}+e^{(\eta-k)t}\big[\rho u\big]_{\text{jump}}-A\big[v^{\gamma+1}\big]_{\text{jump}}e^{(\eta-k)t}e^{kt(\gamma+1)}.
         \label{n=kRanHugRel_3}
         \end{align}
We then substitute (\ref{n=kRanHugRel_3}) into (\ref{n=kRanHugRel_2}) to obtain 
\begin{align}
     \frac{d}{dt}\bigg(\omega_1\omega_\delta\bigg)=& -w_{\delta}e^{(\eta-k)t}\big[\rho u\big]_{\text{jump}}+e^{(\eta-k)t}\big[\rho u^2\big]_{\text{jump}}-A \big[u v^{\gamma+1}\big]_{\text{jump}}e^{(\eta-k)t}e^{kt(\gamma+1)}.\label{n=kRanHugRel_4}
\end{align}
Integration of (\ref{n=kRanHugRel_3})-(\ref{n=kRanHugRel_4}) yields
\begin{equation}
    \label{integrals}
    \begin{cases}
\omega_1\big(t\big)=& -\big[\rho\big]_{\text{jump}}\displaystyle\int_0^t w_{\delta}(s)e^{(\eta-k)s} \ ds+\big[\rho u\big]_{\text{jump}}\displaystyle\int_0^t e^{(\eta-k)s} \ ds-A\big[v^{\gamma+1}\big]_{\text{jump}}\displaystyle\int_0^t e^{(\eta-k)s}e^{ks(\gamma+1)} \ ds,\\ \\
\omega_1\omega_\delta=& -\big[\rho u\big]_{\text{jump}}\displaystyle\int_0^t w_{\delta}(s)e^{(\eta-k)s}\ ds+\big[\rho u^2\big]_{\text{jump}}\displaystyle\int_0^t e^{(\eta-k)s}\ ds-A \big[u v^{\gamma+1}\big]_{\text{jump}}\displaystyle\int_0^t e^{(\eta-k)s}e^{ks(\gamma+1)} \ ds.
    \end{cases}
    \end{equation}
We multiply the first equation with $w_{\delta},$ subtract it from the second and let 
\begin{equation}\label{subs}
    g\big(t\big)=\int_0^t w_{\delta}\big(s\big)e^{(\eta-k)s} \ ds
\end{equation} 
to determine 
\begin{equation}\label{Rankine_n_neqk8}\begin{cases}
&\big(\rho_L-\rho_R\big)g'\big(t\big)g\big(t\big)+g'\big(t\big)\bigg(A\big(\rho_L^{\gamma+1}-\rho_R^{\gamma+1}\big)\bigg(\frac{e^{(k(\gamma+1)+\eta-k)t}-1}{k(\gamma+1)+\eta-k}\bigg)-\big(\rho_L u_L-\rho_Ru_R\big)\bigg(\frac{e^{(\eta-k)t}-1}{\eta-k}\bigg)\bigg)\\&-g\big(t\big)\big(\rho_L u_L-\rho_Ru_R\big)e^{(\eta-k)t}\\&+e^{(\eta-k)t}\bigg(\big(\rho_L u_L^2-\rho_Ru_R^2\big)\bigg(\frac{e^{(\eta-k)t}-1}{\eta-k}\bigg)-\big(Au_L \rho_L^{\gamma+1}-Au_R\rho_R^{\gamma+1}\big)\bigg(\frac{e^{(k(\gamma+1)+\eta-k)t}-1}{k(\gamma+1)+\eta-k}\bigg)\bigg)\\&=0  \hfill \text{when} \ \eta\neq -k\gamma,\ \\ \\
&\big(\rho_L-\rho_R\big)g'\big(t\big)g\big(t\big)+g'\big(t\big)\bigg(A\big(\rho_L^{\gamma+1}-\rho_R^{\gamma+1}\big)t-\big(\rho_L u_L -\rho_Ru_R\big)\bigg(\frac{e^{(\eta-k)t}-1}{\eta-k}\bigg)\bigg)\\&-g\big(t\big)\big(\rho_L u_L-\rho_Ru_R\big)e^{(\eta-k)t}\\&+e^{(\eta-k)t}\bigg(\big(\rho_L u_L^2-\rho_Ru_R^2\big)\bigg(\frac{e^{(\eta-k)t}-1}{\eta-k}\bigg)-\big(Au_L\rho_L^{\gamma+1}-Au_R\rho_R^{\gamma+1}\big)t\bigg)\\&=0 \ \hfill \ \text{when} \ \ \eta= -k\gamma.
\end{cases}\end{equation}
\begin{itemize}
\item When $\eta \neq -\gamma k,$ $\gamma\neq-1$, and $\rho_L=\rho_R$,  the solution of the first ODE in (\ref{Rankine_n_neqk8}) is
\begin{equation}\label{DeltaNNeqK2}\begin{split}\omega_\delta=&\dfrac{-2A\rho_L^\gamma\big(k-\eta\big)e^{t((\gamma+2)k-\eta)}+\big(k\gamma+\eta\big)\big(u_L+u_R\big)e^{2(k-\eta)t}-2\big(u_L+u_R\big)e^{(k-\eta)t}+\frac12\big(u_L+u_R\big)}{2\big(e^{(k-\eta)t}-1\big)^2}\\
&+\dfrac{-2A\big(k-\eta\big)^2\rho_L^\gamma e^{2(k-\eta)t}+A\big(k-\eta\big)\rho_L^\gamma e^{kt(\gamma+1)}}{2\big(\big(\gamma-1\big)k+2\eta\big)\big(e^{(\eta-k)t}-1\big)^2}\end{split}\end{equation}
\item When $\eta = -\gamma k,$ $\gamma\neq-1$, and $\rho_L=\rho_R$, the solution of the second ODE in (\ref{Rankine_n_neqk8}) is
\begin{equation}\label{DeltaNNeqK1}
   \omega_\delta=\frac{1}{\big(e^{(k-\eta)t}-1\big)^2}\bigg(\bigg(\big(1+\big(\eta-k\big)t\big)A\rho_L^\gamma-u_L-u_R\bigg)e^{-2(\eta-k)t}-\big(A\rho_L^{\gamma}+u_L+u_R\big)e^{(k-\eta)t}+u_L+u_R\bigg)
\end{equation}
\end{itemize}
We note that when $\rho_L\neq\rho_R$ the ODEs cannot be solved explicitly. Figure \ref{fig:Delta Shock} is an example obtained numerically when $\eta=-k\gamma$ with parameters $A=-10, \gamma=-4, k=1, \eta=4,\beta=2,\rho_L=2,u_L=3,u_R=2,\rho_R=4$. A similar graph (although it might be flipped across the $t$-axis) is obtained for other parameters as well as for when $\eta\neq -k\gamma.$


\subsubsection{Overcompressible Region}
\begin{figure}[H]
    \centering
    \includegraphics[width=0.7\linewidth]{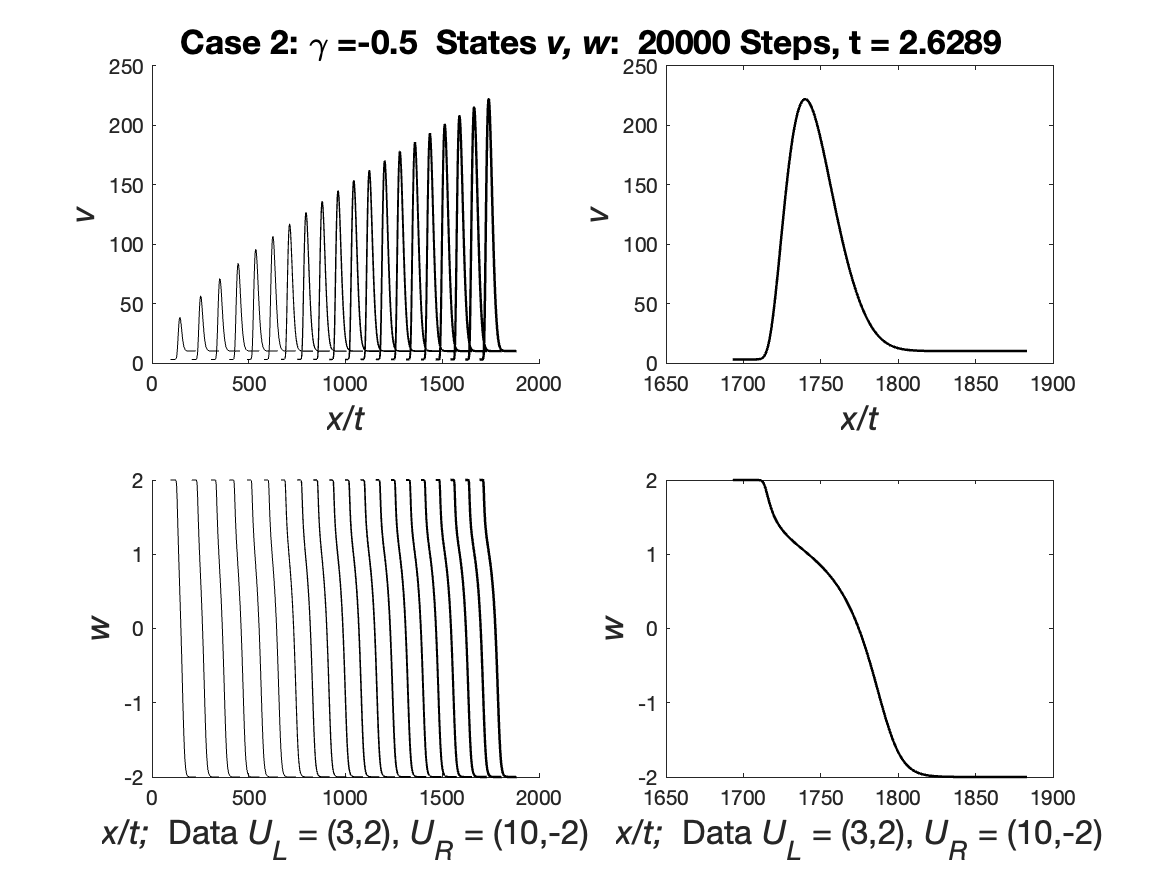}
    \caption{\textit{Region V in Figure \ref{fig:Regions 2}. Dirac delta function in $v$. Parameters: $\gamma = -0.5, A = -10, \eta = 3, k = -0.01, \beta = 10$.}}
    \label{fig:Case2Delta}
\end{figure}

We seek delta-shocks connecting a given left state $(v_-,w_+)$ with a right state $(v,w)$ that are overcompressive, meaning that all characteristic curves run into the delta-shock curve from both sides. Therefore, we require the following inequality 
\begin{equation}\label{OCompIneq}
    \lambda_1\big(v,w\big) < \lambda_2\big(v,w\big) \leq \frac{dx\big(t\big)}{dt} \leq \lambda_1\big(v_-,w_-\big) < \lambda_2\big(v_-,w_-\big). 
\end{equation}
The outer inequalities always hold. Note, $\omega_\delta$ generally cannot be solved explicitly. Thus, we consider $\lambda_2(v,w)<\lambda_1(v_-,w_-)$ to locate the region which would contain the right states that result in a strictly overcompressive delta-shock. This inequality indicates an upper border \begin{equation}\label{OCompBound}
   J: w = w_{-} -A \big(\gamma +1\big)v_-^{\gamma} e^{kt(\gamma +1)} + A v^{\gamma} e^{kt(\gamma +1)}.
\end{equation}
After some simplification, this curve is the same in the $\eta=k$ case. Hence, it yields the same region. $J$ is above $S_{\delta},$ given by (\ref{n=kbound}), when $-1<\gamma<0$ and below when $\gamma<-1$. As mentioned above, Dirac delta functions are observed numerically only in $v$, as shown in Figure \ref{fig:Case2Delta}.

\subsubsection{Region Shift for $k(\gamma+1)<0$}
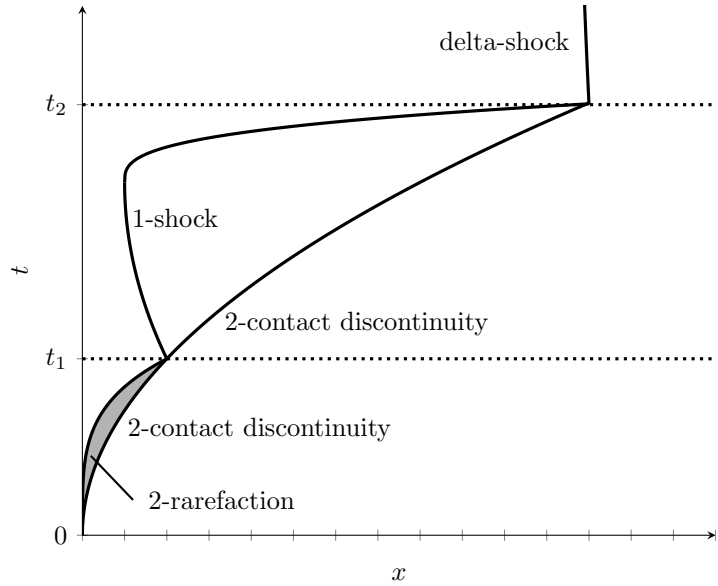
\begin{figure}[H]
    \centering
    \begin{tikzpicture}[x=100,y=50]
    \begin{axis}[
      ymax=3,
      ymin=0,
      xmax=15,
      xmin=0,
      axis x line=bottom,
      axis y line=left,
      domain=0:10,
      xtick={0, 1, ..., 15},
      ytick={0, 1, 2.44, ..., 3},
      samples=1001,
      yticklabels={$0$,$t_1$,$t_2$,},
      xticklabels=\empty,
      xlabel = \(x\),
      ylabel = {\(t\)}
    ]
    \addplot[domain=0:2,
    name path = A,
    very thick
    ] {x^(1/4)/2^(1/4)};
    \addplot[domain=0:2,
    name path = B,
    very thick
    ] {x^(1/2)/(2^(1/2)};
    \draw[color=black,thick] (axis cs:.2,.45) -- (axis cs:1.2,.2);
    \addplot [black!30] fill between [of = A and B, soft clip={domain=-5:0}];
    \node at (axis cs:3.3,.2)(B){\color{black}2-rarefaction};
    \node at (axis cs:4.2,.6)(B){\color{black}2-contact discontinuity};
    \node at (axis cs:6.5,1.2)(B){\color{black}2-contact discontinuity};
    \addplot[domain=0:15,
    dotted,
    very thick
    ] {1};
    \addplot[domain=2:12,
    very thick
    ] {x^(1/2)/(2^(1/2)};
    \addplot[domain=10:12,
    very thick
    ] {2.44*(e^(-2*x)/e^(-2*12))};
    \addplot[domain=0:15,
    dotted,
    very thick
    ] {2.44};
    \addplot[domain=1:12,
    very thick] {(x-1)^(1/3)/5+2};
    \addplot[domain=1:2,
    very thick] {-(x-1)^(1/2)+2};
    \node at (axis cs:2.2,1.8)(B){\color{black}1-shock};
    \node at (axis cs:10,2.8)(B){\color{black}delta-shock};
  \end{axis}
\end{tikzpicture}
    \caption{An example solution for a $(V_R,W_R)$ originally in Region II Case II}
    \label{fig:SolExamp}
\end{figure}
\begin{figure}[H]
\begin{center}
\begin{tikzpicture}
  \begin{axis}[
      ymax=5,
      ymin=-5,
      xmax=5,
      xmin=0,
      axis x line=middle,
      axis y line=left,
      domain=0:10,
      xtick={0, 1, ..., 5},
      ytick={0, 1, ..., 5},
      samples=1001,
      yticklabels=\empty,
      xticklabels=\empty,
      xlabel = \(v\),
      ylabel = {\(w\)}
    ]
    \addplot[domain=0:5,
    color=black,
    thick
    ] {3};
    \draw[color=black,thick] (axis cs:.05,-5) -- (axis cs:.05,3);
    \draw[color=black] (axis cs:.025,-1) -- (axis cs:.5,-.5);
    \node at (axis cs:0.75,-0.5)(B){\color{black}VII};
    \node at (axis cs:2.5,3.4)(B){\color{black}VI};
    \node at (axis cs:2.5,2.7)(B){\color{black}V};
  \end{axis}
\end{tikzpicture}
\end{center}
\caption{\textit{Limit behavior of the regions for $\gamma<-1$ and $k>0$}}
    \label{fig:Time Regions 1}
\end{figure}
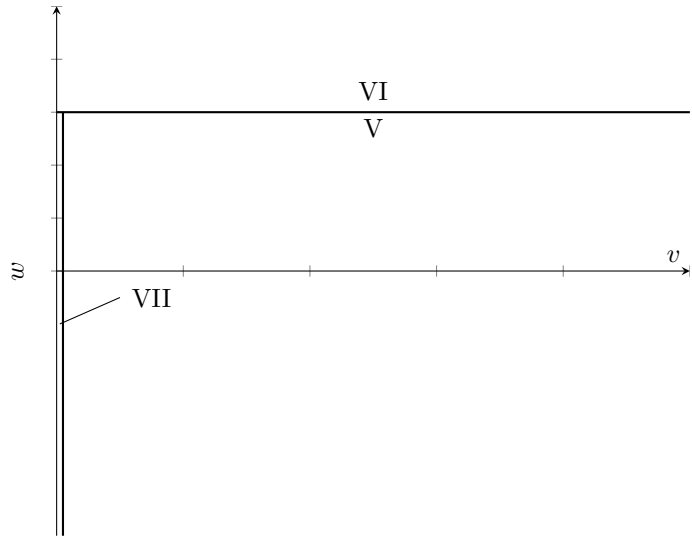

In this case, as $t\to\infty$, the $C_2$ and $R_2$ curves converge to $w=w_L$ with $C_2$ maintaining an asymptote at $v=0$. This allows for a point in a given region to shift to another as time progresses. The behavior of regions is found to progress as shown in Figure \ref{fig:Time Regions 1} and Figure \ref{fig:Time Regions 2}.
In addition, as time approaches infinity, (\ref{n=kbound}) and (\ref{OCompBound}) will approach $w = w_L$. Therefore the set of overcompressible points will be a subset of Region V in Figure \ref{fig:Time Regions 1}. The solutions to the Riemann problem will consist of:
\begin{itemize}
\item A 1-shock followed by a 2-contact discontinuity when the right state is in Region VI.
\item A 1-shock followed by a 2-rarefaction when the right state is in Region VII.
\item An overcompressive delta-shock or a combination of a delta-shock and a classical wave when the right state is in Region V.
\end{itemize}
The limit will affect Case 2 in a similar way:
\begin{figure}[H]
\begin{center}
\begin{tikzpicture}
  \begin{axis}[
      ymax=5,
      ymin=-5,
      xmax=5,
      xmin=0,
      axis x line=middle,
      axis y line=left,
      domain=0:10,
      xtick={0, 1, ..., 5},
      ytick={0, 1, ..., 5},
      samples=1001,
      yticklabels=\empty,
      xticklabels=\empty,
      xlabel = \(v\),
      ylabel = {\(w\)}
    ]
    \addplot[domain=.05:5,
    color=black,
    thick
    ] {3};
    \draw[color=black,thick] (axis cs:.05,-5) -- (axis cs:.05,3);
    \node at (axis cs:2.5,3.4)(B){\color{black}I};
    \node at (axis cs:2.5,2.7)(B){\color{black}V};
  \end{axis}
\end{tikzpicture}
\end{center}
\caption{\textit{Limit behavior of the regions for $-1<\gamma<0$ and $k<0$}}
    \label{fig:Time Regions 2}
\end{figure}
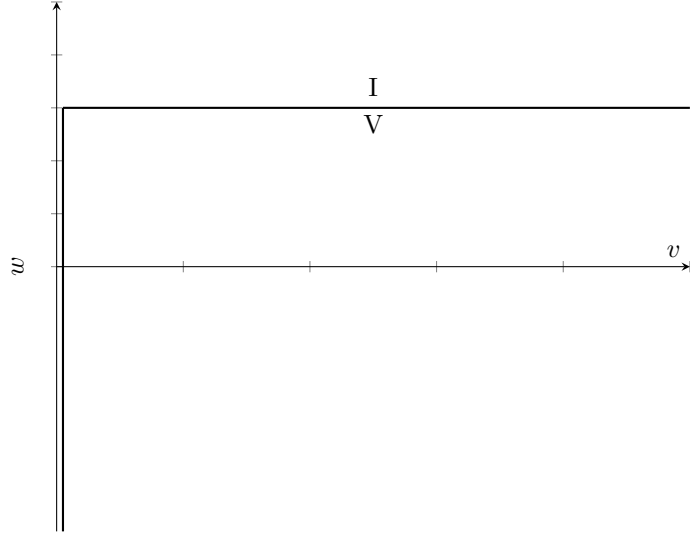
However, we now expect: 
\begin{itemize}
\item A 2-rarefaction followed by a 2-contact discontinuity when the right state is in Region I.
\item An overcompressive delta-shock or a combination of a delta-shock and a classical wave when the right state is in Region V. 
\end{itemize}

\begin{figure}[H]
\centering
\includegraphics[width=0.7\linewidth]{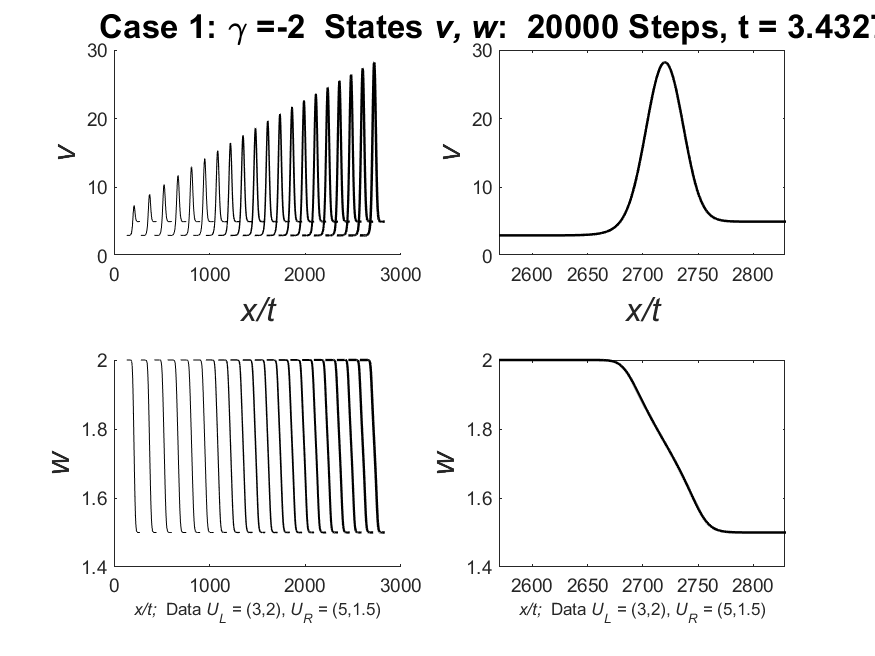}
\caption{\textit{Region V of Figure \ref{fig:Time Regions 1}. The region shifts to a delta-shock over time in $v$. Parameters: $\gamma = -2, A = -10, \eta = 3, k = 2, \beta = 10$.}}
\label{fig:Case1OverComp}
\end{figure}

\begin{figure}[H]
\centering
\includegraphics[width=0.7\linewidth]{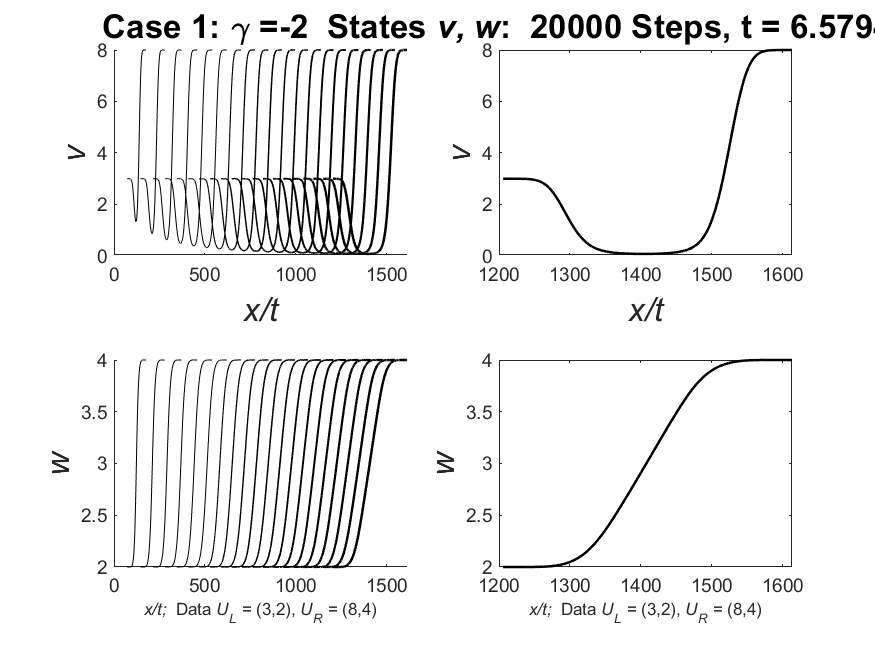}
\caption{\textit{Region VI of Figure \ref{fig:Time Regions 1}. Region shifts to $S_1C_2$. Parameters: $\gamma = -2, A = -10, \eta = 3, k = 2, \beta = 10$.}}
\label{fig:Case1_RegionVI_Collapse}
\end{figure}

\begin{figure}[H]
\centering
\includegraphics[width=0.7\linewidth]{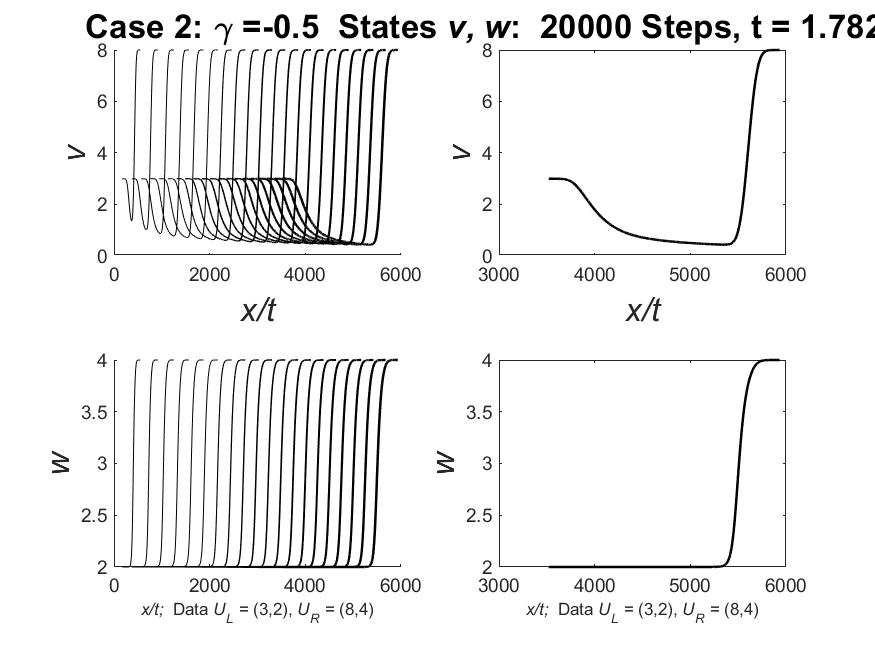}
\caption{\textit{Region I of \ref{fig:Time Regions 2}. Region shifts to $R_2C_2$. Parameters: $\gamma = -0.5, A = -10, \eta = 3, k = -2, \beta = 10$.}}
\label{fig:Case2_RegionVI_Collapse}
\end{figure}

\subsubsection{Region Shift for $k(\gamma+1)>0$}
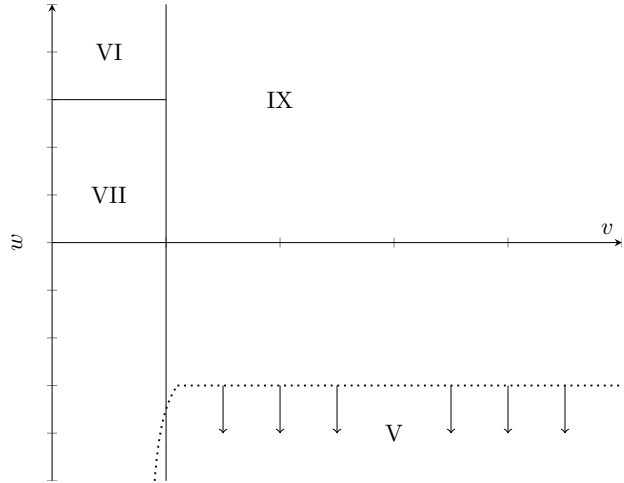
\begin{figure}[H]
\begin{center}
\begin{tikzpicture}[scale=0.9]
  \begin{axis}[
      axis x line=middle,
      axis y line=left,
      ymax=5,
      ymin=-5,
      xmax=5,
      xmin=0,
      domain=0:10,
      xtick={0, 1, ..., 5},
      ytick={-5, -4, ..., 5},
      samples=1001,
      yticklabels=\empty,
      xticklabels=\empty,
      xlabel = \(v\),
      ylabel = {\(w\)}
    ]
    \addplot [
    domain=0:1,
    color=black,
    ] {3};
    \addplot[
    domain=1.1:5,
    color = black,
    dotted,
    thick,
    ] {-3};
    \addplot[
    domain=.9:1.1,
    color = black,
    dotted,
    thick,
    ] {-1/((x-.8)/.3)-2};
    \draw[color=black] (axis cs:1,-5) -- (axis cs:1,5);
    \node at (axis cs:.5,1)(B){\color{black}VII};
    \node at (axis cs:.5,4)(B){\color{black}VI};
    \node at (axis cs:2,3)(B){\color{black}IX};
    \node at (axis cs:3,-4)(B){\color{black}V};
    \draw[to-] (axis cs: 1.5,-4) -> (axis cs:1.5,-3);
    \draw[to-] (axis cs: 2,-4) -> (axis cs:2,-3);
    \draw[to-] (axis cs: 2.5,-4) -> (axis cs:2.5,-3);
    \draw[to-] (axis cs: 3.5,-4) -> (axis cs:3.5,-3);
    \draw[to-] (axis cs: 4,-4) -> (axis cs:4,-3);
    \draw[to-] (axis cs: 4.5,-4) -> (axis cs:4.5,-3);
  \end{axis}
\end{tikzpicture}
\end{center}
\caption{\textit{Limit behavior of the regions for $\gamma<-1$ and $k<0$}}
    \label{fig:Time Regions 4}
\end{figure}

As $t\to\infty$ the shift of the regions will require more careful consideration since 
\begin{equation}\label{n=kRegionLim}\lim_{t\to\infty}\bigg(w_{-}+\frac{A}{v}\big(ve^{kt}\big)^{\gamma+1}\bigg)=-\infty.\end{equation}
Therefore, $S_{\delta}$ moves farther down in the plane (when $\gamma<-1$ the overcompressive region moves in the same spirit because it's located below $S_{\delta}$). Next, we check the time behavior of $C_2.$ Since 
\begin{equation}\label{n=kRegionLim2}\lim_{t\to\infty}\bigg(w_-+\frac{A}{v}\big(ve^{kt}\big)^{\gamma+1}-\frac{A}{v_-}\big(v_-e^{kt}\big)^{\gamma+1}\bigg)=\lim_{t\to\infty}\bigg(w_-+A\big(v^{\gamma}-v_-^{\gamma}\big)e^{kt(\gamma+1)}\bigg)=\begin{cases}
-\infty \ \ \text{when} \ \ v<v_-,\\
+\infty \ \ \text{when} \ \ v>v_-,
    \end{cases}
\end{equation}
the regions will shift as shown in Figures \ref{fig:Time Regions 3} and \ref{fig:Time Regions 4}.

\vspace{5mm}

When $-1<\gamma<0$ the solutions to the Riemann problem will consist of:
\begin{itemize}
    \item A 2-rarefaction followed by a 2-shock when the right state is in Region I. See Figure \ref{fig:Case4_IV_SC}.
    \item A 1-shock followed by a 2-shock when the right state is in Regions III and IV. See Figure \ref{fig:Case4Region2Collapse}.
\end{itemize}
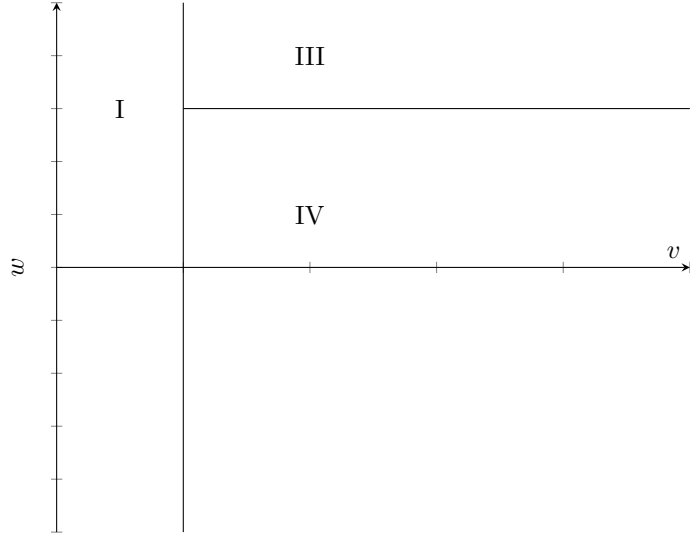
\begin{figure}[H]
\begin{center}
\begin{tikzpicture}
  \begin{axis}[
      axis x line=middle,
      axis y line=left,
      ymax=5,
      ymin=-5,
      xmax=5,
      xmin=0,
      domain=0:10,
      xtick={0, 1, ..., 5},
      ytick={-5, -4, ..., 5},
      samples=1001,
      yticklabels=\empty,
      xticklabels=\empty,
      xlabel = \(v\),
      ylabel = {\(w\)}
    ]
    \addplot [
    domain=1:5,
    color=black,
    ] {3};
    \draw[color=black] (axis cs:1,-5) -- (axis cs:1,5);
    \node at (axis cs:2,1)(B){\color{black}IV};
    \node at (axis cs:2,4)(B){\color{black}III};
    \node at (axis cs:.5,3)(B){\color{black}I};
  \end{axis}
\end{tikzpicture}
\end{center}
\caption{\textit{Limit behavior of the regions for $-1<\gamma<0$ and $k>0$}}
    \label{fig:Time Regions 3}
\end{figure}


\begin{figure}[H]
\centering
\includegraphics[width=0.7\linewidth]{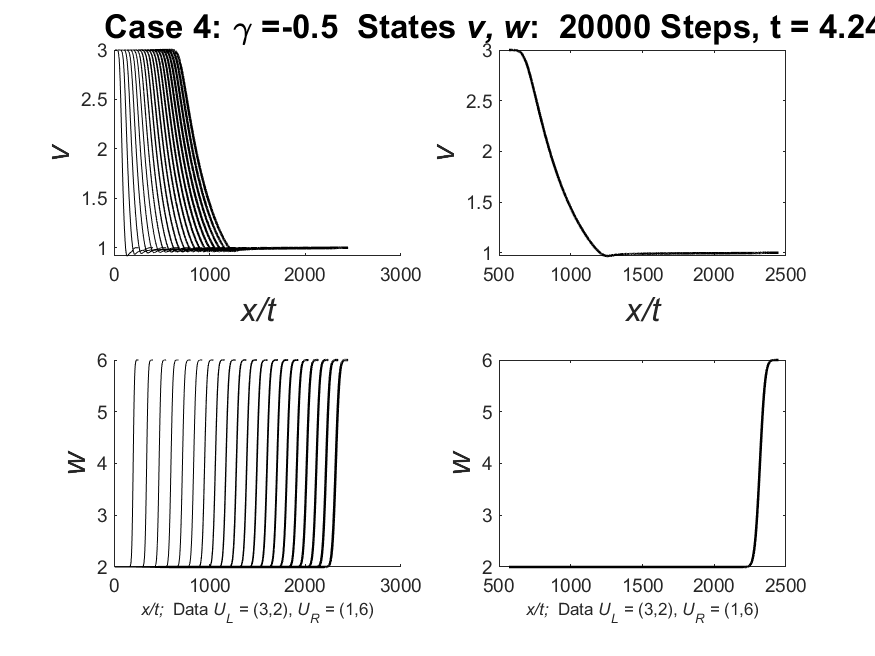}
\caption{\textit{Region I of Figure \ref{fig:Time Regions 3}. Region shifts to $R_2C_2$ in time. Parameters: $\gamma = -0.5, A = -10, \eta = 3, k = 0.6, \beta = 10$.}}
\label{fig:change later}
\end{figure}

\begin{figure}[H]
\centering
\includegraphics[width=0.7\linewidth]{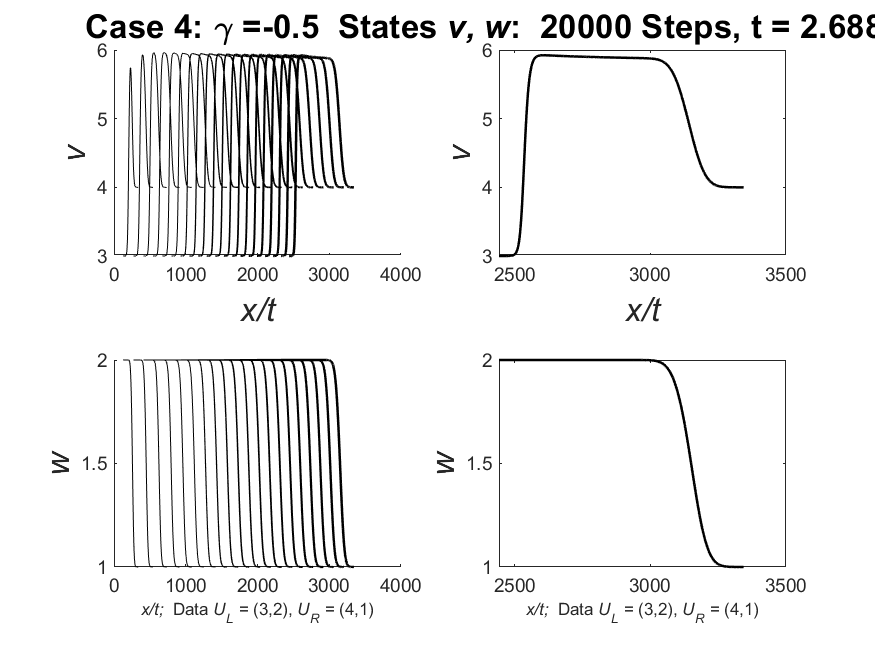}
\caption{\textit{Region III of Figure \ref{fig:Time Regions 3}. Region shifts to $S_1C_2$ in time. Parameters: $\gamma = -0.5, A = -10, \eta = 3, k = 0.6, \beta = 10$.}}
\label{fig:Case4Region2Collapse}
\end{figure}
The $\gamma<-1$ case has a similar transformation as seen below.
In this final case, the solution will consist of:
\begin{itemize}
    \item A 1-shock followed by a 2-rarefaction when $(v_R,w_R)$ is in Region VI. This means that $R_2$ remains above $C_2$ during the region shift. 
    \item A 1-shock followed by a 2-contact discontinuity when $(v_R,w_R)$ is in Region VII. 
    \item Either a delta-shock followed by a 2-wave or a 2-contact discontinuity followed by a delta-shock when the right state is in Region IX. The set of overcompressible points will be a subset of Region V in Figure \ref{fig:Time Regions 4}.
\end{itemize}

\vspace{-0.8cm}
\begin{figure}[H]
\centering
\includegraphics[width=0.7\linewidth]{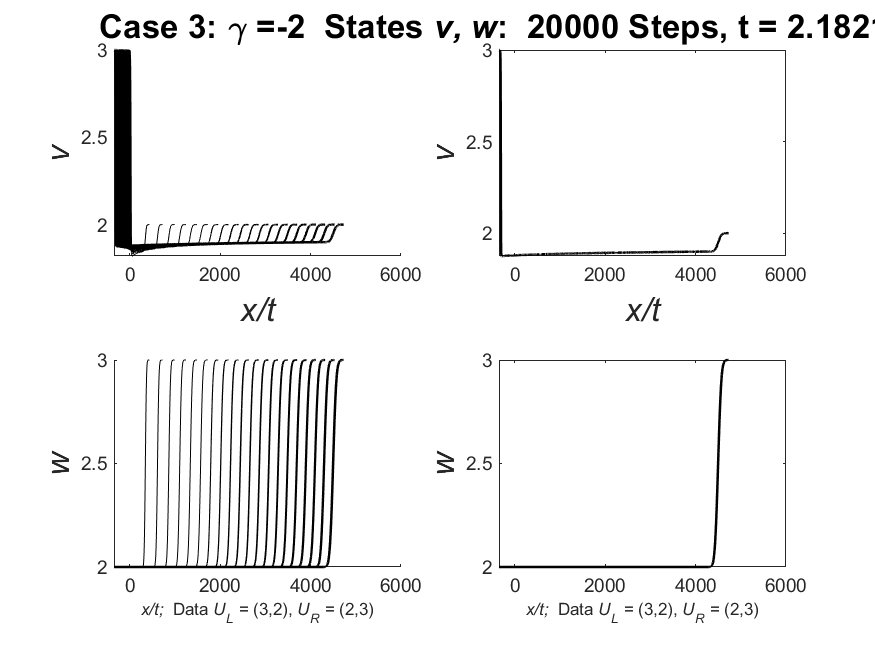}
\caption{\textit{Region VI of Figure \ref{fig:Time Regions 4}. Region shifts to $S_1R_2$. Parameters: $\gamma = -2, A = -10, \eta = 3, k = -0.6, \beta = 10$.}}
\label{fig:Case3Region6Collapse}
\end{figure}
\vspace{-0.8cm}
\begin{figure}[H]
\centering
\includegraphics[width=0.7\linewidth]{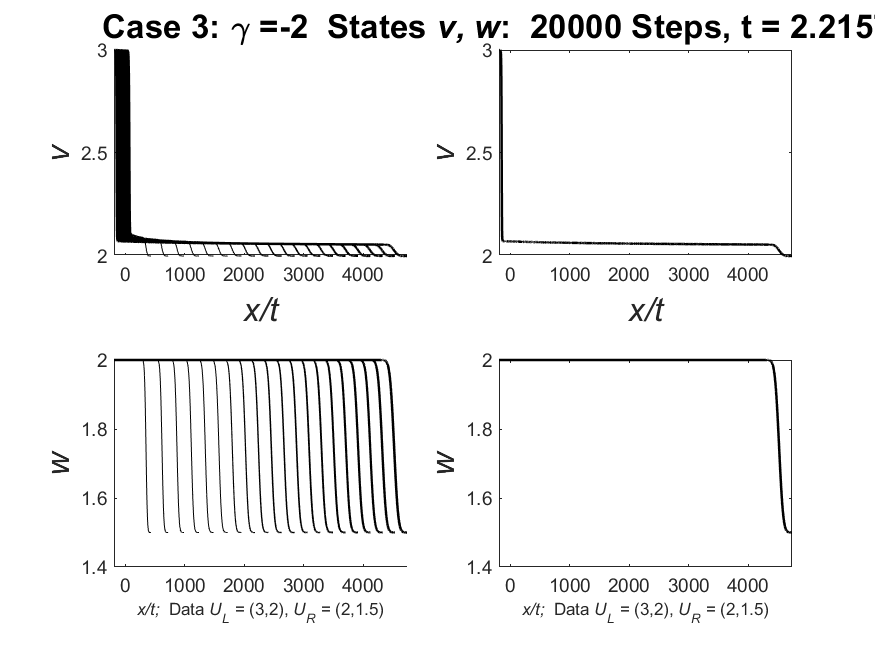}
\caption{\textit{Region VII in Figure \ref{fig:Time Regions 4}. Region shifts to $S_1C_2$. Parameters: $\gamma = -2, A = -10, \eta = 3, k = -0.6, \beta = 10$.}}
\label{fig:Case3Region7Collapse}
\end{figure}
\vspace{-0.8cm}
\subsubsection{Non-overcompressive Regions}

\begin{figure}[H]
\centering
\includegraphics[width=0.7\linewidth]{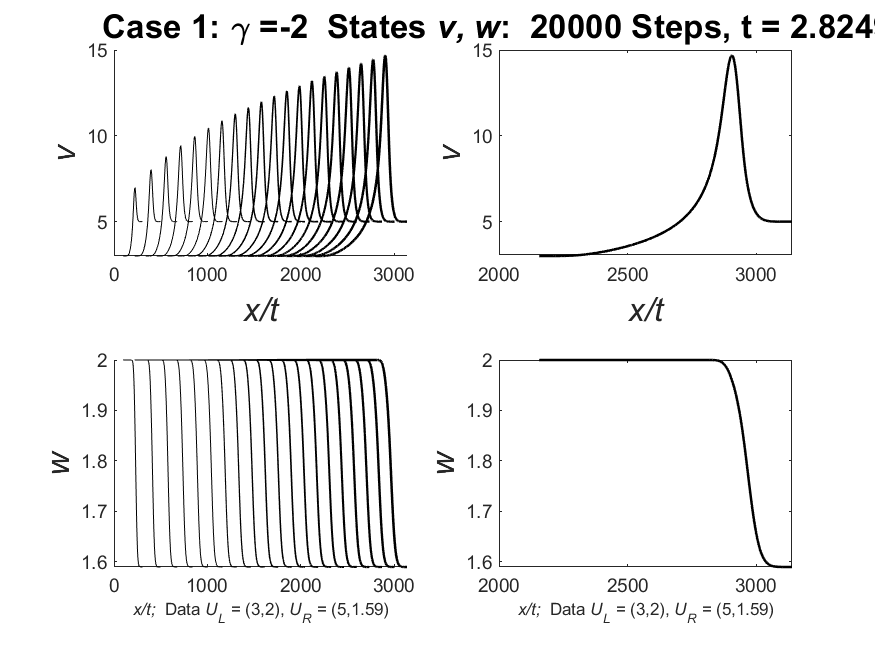}
\caption{\textit{Region IX of Figure \ref{fig:Regions 4}, $\delta 2$-$\text{wave}$, between $S_\delta$ and $S_o$. Parameters: $\gamma = -2, A = -10, \eta = 3, k = 0.01, \beta = 10$.}}
\label{fig:Case1BoundaryComboDS}
\end{figure}

\begin{figure}[H]
\centering
\includegraphics[width=0.7\linewidth]{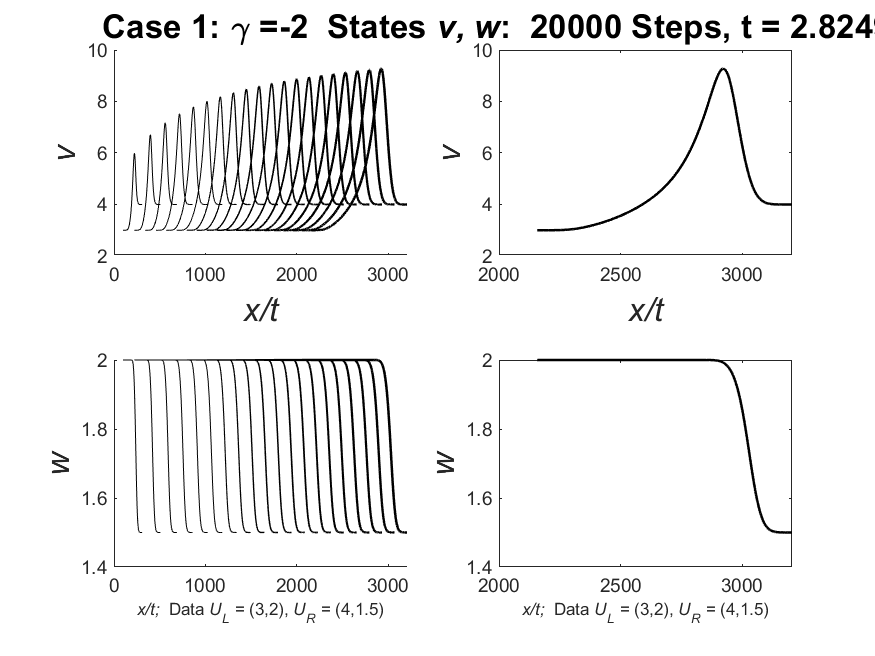}
\caption{\textit{Region IX of Figure \ref{fig:Regions 4}, $C_2\delta$, close to and above $S_\delta$. Parameters: $\gamma = -2, A = -10, \eta = 3, k = 0.01, \beta = 10$.}}
\label{fig:Case1ComboAwayC2}
\end{figure}

\begin{figure}[H]
\centering
\includegraphics[width=0.7\linewidth]{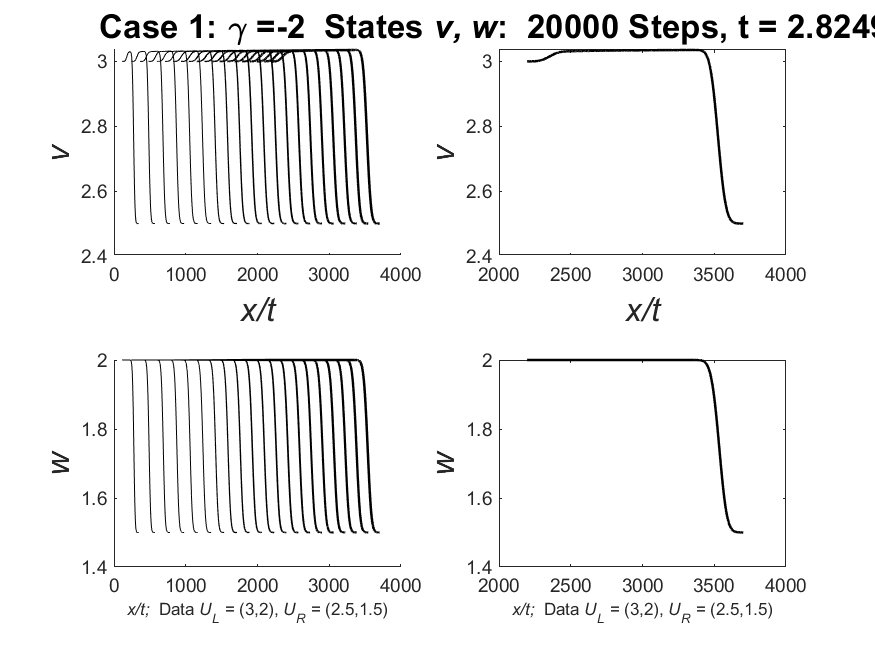}
\caption{\textit{Region IX of Figure \ref{fig:Regions 4}, $C_2\delta$, farther from and above $S_\delta$. Parameters: $\gamma = -2, A = -10, \eta = 3, k = 0.01, \beta = 10$.}}
\label{fig:Case1ComboCloseC2}
\end{figure}

Before region shift, solutions with right states in Region IX, Case 1, display combinations of a delta-shock and a classical solution. Between $S_\delta$ and $S_o$, a delta-shock is followed by a 2-wave. This is supported by Figure \ref{fig:Case1BoundaryComboDS}, as the 2-wave in $w$ only occurs after the initial delta-shock. Above $S_\delta$, we observe a 2-contact discontinuity followed by a delta-shock. Note that the diffusion in the first wave is likely due to the cell averaging involved in the LLF method as the delta-shock grows. Further testing with other numerical methods is needed to confirm this. We find that proximity to $S_\delta$ shows increased delta-shock characteristics numerically, while proximity to $C_2$ muddles those characteristics. We demonstrate this in Figure \ref{fig:Case1ComboAwayC2} and Figure \ref{fig:Case1ComboCloseC2}. We propose this is due to the curve needing to travel farther along $C_2$ to a middle state, raising the overcompressive region, allowing for a delta-shock to be taken. If this is the case, it would explain why numerically a proximity to the $C_2$ curve would cause an unclear delta characteristic as the curve must travel incredibly far in order to raise the overcompressive region a sufficient amount.
\vspace{-0.1cm}
\vspace{-1cm}

\begin{figure}[H]
\centering
\includegraphics[width=0.7\linewidth]{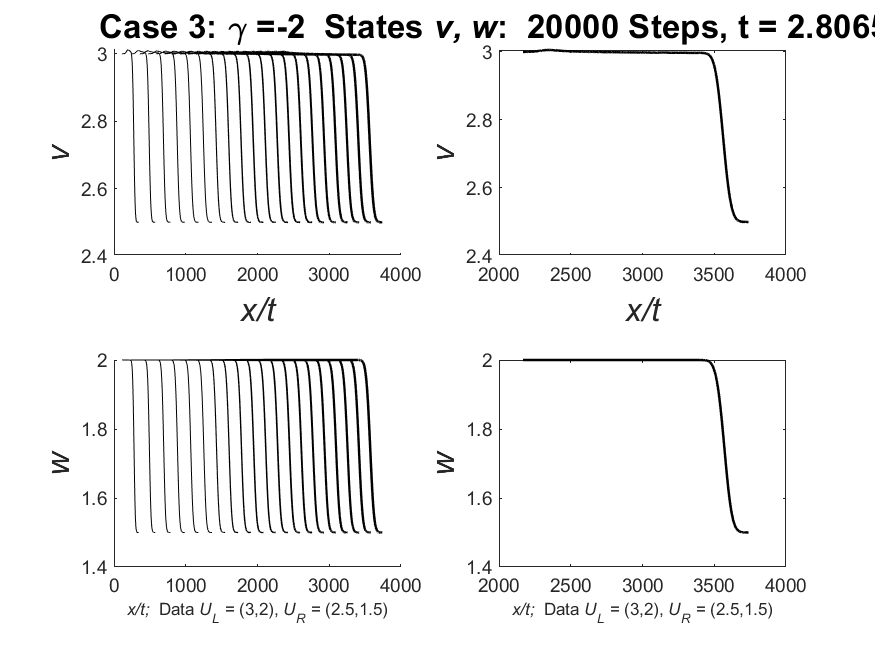}
\caption{\textit{Region IX of Figure \ref{fig:Regions 3}, $C_2\delta$, close to $C_2$ and above $S_\delta$. Parameters: $\gamma = -2, A = -10, \eta = 3, k = -0.01, \beta = 10$.}}
\label{fig:Case3ComboCloseC2Smallk}
\end{figure}

\begin{figure}[H]
\centering
\includegraphics[width=0.7\linewidth]{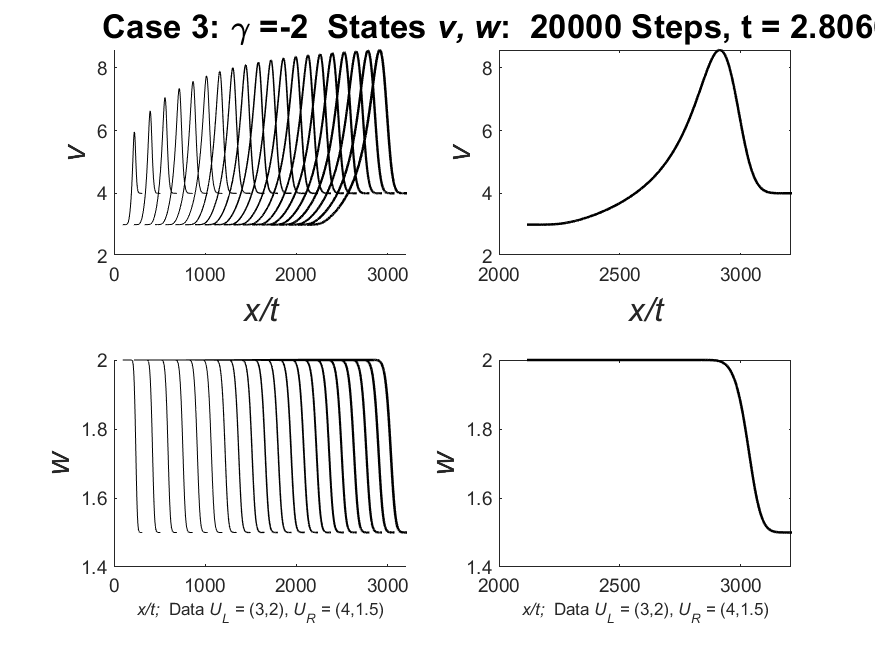}
\caption{\textit{Region IX of Figure \ref{fig:Regions 3}, $C_2\delta$, far from $C_2$ and above $S_\delta$. Parameters: $\gamma = -2, A = -10, \eta = 3, k = -0.01, \beta = 10$.}}
\label{fig:Case3ComboAwayC2smallk}
\end{figure}

\begin{figure}[H]
\centering
\includegraphics[width=0.7\linewidth]{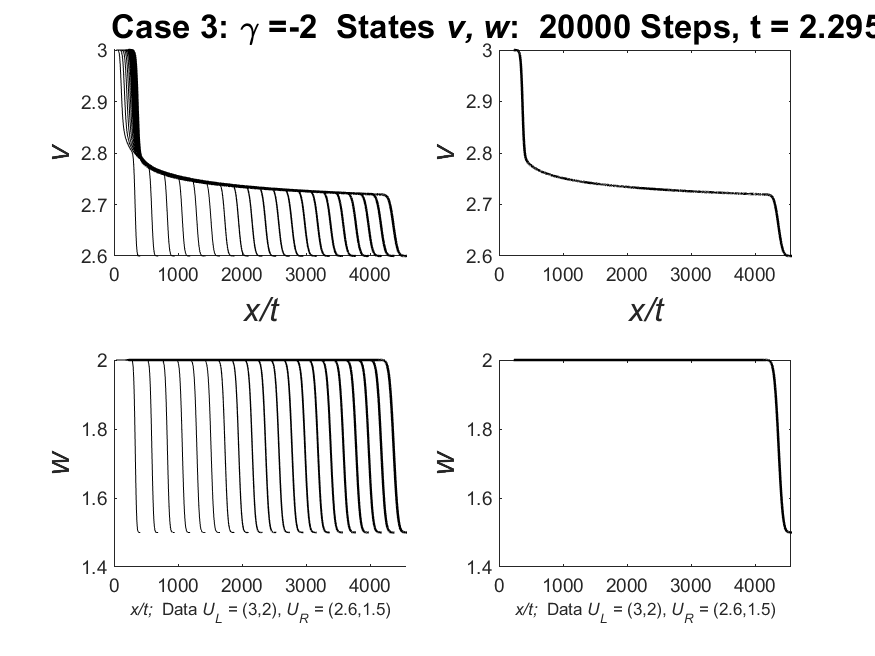}
\caption{\textit{Region IX of Figure \ref{fig:Time Regions 4}, $C_2\delta$, close to $C_2$ and above $S_\delta$. Parameters: $\gamma = -2, A = -10, \eta = 3, k = -2, \beta = 10$.}}
\label{fig:Case3ComboCloseC2largek}
\end{figure}

\begin{figure}[H]
\centering
\includegraphics[width=0.7\linewidth]{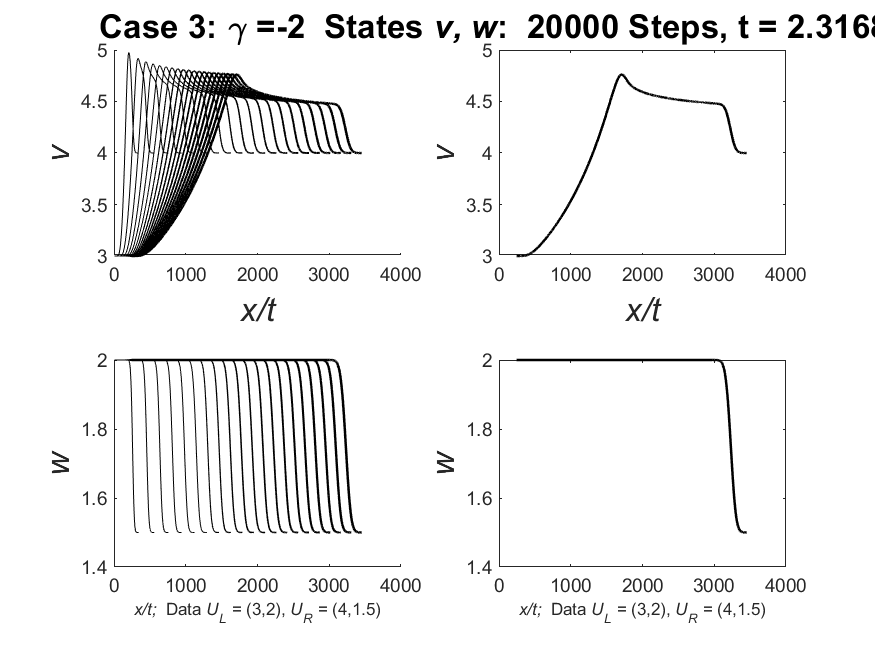}
\caption{\textit{Region IX of \ref{fig:Time Regions 4}, $C_2\delta$, far from $C_2$ and above $S_\delta$. Parameters: $\gamma = -2, A = -10, \eta = 3, k = -2, \beta = 10$.}}
\label{fig:Case3ComboAwayC2largek}
\end{figure}

Note that in Case 3, $S_\delta$ and $S_o$ approach the same limit. Case 1 shares the same short term behavior of non-overcompressive delta regions, but it is lost upon region shift. This is shown in Figures \ref{fig:Case3ComboCloseC2Smallk} \ref{fig:Case3ComboAwayC2smallk}, \ref{fig:Case3ComboCloseC2largek}, and \ref{fig:Case3ComboAwayC2largek}. As shown in Figure \ref{fig:Case3ComboCloseC2largek}, as overcompressibility disappears during region shift, the delta shock becomes much weaker.  

\subsection{$\eta = k$ Case}
We follow the delta-shock definition as before, therefore our solution should satisfy the equations
\begin{equation}\label{n=kDeltaDefP2}
    \begin{split}
        &\bigg\langle v,\phi_t\bigg\rangle+\bigg\langle v\big(w+\beta t\big)-A\big(ve^{kt}\big)^{\gamma+1},\phi_x\bigg\rangle=0\\
&\bigg\langle vw,\phi_t\bigg\rangle+\bigg\langle v\big(w+\beta t\big)w-A\big(ve^{kt}\big)^{\gamma+1}w,\phi_x\bigg\rangle=0
    \end{split}
\end{equation}
for any $\phi \in C_{0}^{\infty}\big(\mathbb{R} \times \mathbb{R}^+\big)$. We use the properties of the Dirac delta function in a similar fashion to the $\eta \neq k$ case to observe that
\begin{equation}\label{n=kGreensT1}
   \begin{split}
      \bigg\langle v,\phi_t\bigg\rangle+\bigg\langle v\big(w+\beta t\big)-A\big(ve^{kt}\big)^{\gamma+1},\phi_x\bigg\rangle
    =   &\int_0^\infty \int_{-\infty}^{x(t)}v_L\phi_tdxdt+  \int_0^\infty\int_{-\infty}^{x(t)} v_L\big(w_L + \beta t\big)\phi_x dxdt \\&+\int_0^\infty \int_{x(t)}^\infty v_R\big(w_R + \beta t\big)\phi_x dxdt+\int_0^\infty \int_{x(t)}^\infty v_R\phi_tdxdt \\&+ \int_0^\infty\bigg(\omega_1 \phi_t + \omega_1\big(w_{\delta} + \beta t\big) \phi_x\bigg) dt \\ & -\int_0^\infty\int_{-\infty}^{x(t)} A \big(v_L e^{kt}\big)^{\gamma +1} \phi_x dxdt -\int_0^\infty\int_{x(t)}^{\infty} A \big(v_R e^{kt}\big)^{\gamma +1} \phi_x dxdt 
    \end{split}
    \end{equation}
and 
\begin{equation}\label{n=kGreensT2}
 \begin{split}
       \bigg\langle vw,\phi_t\bigg\rangle+\bigg\langle v(w+\beta t)w-A(ve^{kt})^{\gamma+1}w,\phi_x\bigg\rangle
    =   &\int_0^\infty \int_{-\infty}^{x(t)}v_Lw_L\phi_tdxdt+  \int_0^\infty\int_{-\infty}^{x(t)} v_L(w_L + \beta t)w_L\phi_x dxdt \\&+\int_0^\infty \int_{x(t)}^\infty v_R(w_R + \beta t)w_R\phi_x dxdt+\int_0^\infty \int_{x(t)}^\infty v_Rw_R\phi_tdxdt \\&+ \int_0^\infty\bigg(\omega_1 w_{\delta}\phi_t + \omega_1(w_{\delta} + \beta t)w_{\delta} \phi_x\bigg) dt \\ & -\int_0^\infty\int_{-\infty}^{x(t)} A \big(v_L e^{kt}\big)^{\gamma +1}w_L \phi_x dxdt \\ &-\int_0^\infty\int_{x(t)}^{\infty} A \big(v_R e^{kt}\big)^{\gamma +1}w_R \phi_x dxdt.
    \end{split}
    \end{equation}
We require 
\begin{equation}
    \label{shock_2}
\frac{dx\big(t\big)}{dt} = \sigma\big(t\big) = \omega_\delta\big(t\big)+\beta t
\end{equation} 
and use Green's theorem to get
\begin{equation}\label{n=kGreensT21}
    \begin{split}
        \bigg\langle v,\phi_t\bigg\rangle+\bigg\langle v\big(w+\beta t\big)-A\big(ve^{kt}\big)^{\gamma+1},\phi_x\bigg\rangle
        = -&\oint v_L\phi dx +\bigg(v_L\big(w_L + \beta t\big) -A\big(v_L e^{kt}\big)^{\gamma +1}\bigg) \phi dt \\ -&\oint -v_R\phi dx +\bigg(v_R\big(w_R + \beta t\big) -A\big(v_R e^{kt}\big)^{\gamma +1}\bigg) \phi dt\\ + &\int_0^\infty \omega_1 d\phi =0
    \end{split}
\end{equation}
and 
\begin{equation}\label{n=kGreensT22}
    \begin{split}
        \bigg\langle vw,\phi_t\bigg\rangle+\bigg\langle v\big(w+\beta t\big)w-A\big(ve^{kt}\big)^{\gamma+1}w,\phi_x\bigg\rangle
        = -&\oint v_Lw_L\phi dx +\bigg(v_L\big(w_L + \beta t\big)w_L -A\big(v_L e^{kt}\big)^{\gamma +1}w_L\bigg) \phi dt \\ -&\oint -v_Rw_R\phi dx +\bigg(v_R\big(w_R + \beta t\big)w_R -A\big(v_R e^{kt}\big)^{\gamma +1}w_R\bigg) \phi dt\\ + &\int_0^\infty \omega_1w_{\delta} d\phi = 0
    \end{split}
\end{equation}
If we also require
\begin{equation}\label{n=kRanHugRel}\begin{cases}\dfrac{d\omega_1}{dt}= \bigg[v\big(w + \beta t\big) -A\big(v e^{kt}\big)^{\gamma +1}\bigg]_{\text{jump}} -\big[v\big]_{\text{jump}}\sigma\big(t\big),
\\ \\\dfrac{d}{dt}\big(\omega_1\omega_\delta\big)= \bigg[vw\big(w + \beta t\big) -Aw\big(v e^{kt}\big)^{\gamma +1}\bigg]_{\text{jump}} -\big[vw\big]_{\text{jump}}\sigma\big(t\big)\end{cases}\end{equation}
then $(v,w)$ satisfies the system in the sense of distributions. Similar to the $\eta \neq k$ case, we return to the original variables $\rho_L$, $\rho_R$, $u_L$ and $u_R,$ integrate the equations, subtract one from the product of the other, and let $g(t)=\int_0^t \omega_\delta(s) \ ds$ to form the ODE:
\begin{equation}\label{n=kODE}
    \begin{split}
        &-\big(\rho_L -\rho_R\big) g'\big(t\big)g\big(t\big)+g'\big(t\big)\bigg(\big(\rho_Lu_L - \rho_Ru_R\big)t - A\bigg(\frac{e^{k(\gamma +1)t}-1}{k(\gamma+1)}\bigg) \big(\rho_L^{\gamma +1} -\rho_R^{\gamma +1}\big)\bigg)\\ &+\big(\rho_Lu_L - \rho_Ru_R\big)g\big(t\big)-\big(\rho_Lu_L^2 - \rho_Ru_R^2\big)t+A\bigg(\frac{e^{k(\gamma +1)t}-1}{k(\gamma +1)}\bigg) \big(\rho_L^{\gamma +1}u_L -\rho_R^{\gamma +1}u_R\big) = 0.
    \end{split}
\end{equation}
When $\rho_L= \rho_R$, this equation can be solved explicitly:
\begin{equation}\label{n=kODESol}
    \omega_\delta = \frac{1}{2}\big(u_L +u_R\big) - \frac{kt(\gamma+1)\rho_L^{\gamma}Ae^{k(\gamma +1)t}-\rho_L^{\gamma}Ae^{k(\gamma +1)t}+A\rho_L^\gamma}{t^2 k^2(\gamma +1)^2 }
\end{equation}
When $\gamma\neq -1$ and $\rho_L\neq\rho_R$ the equation is solved numerically to be similar to Figure \ref{fig:Delta Shock}. Furthermore, the overcompressible regions are based on a simplified version of (\ref{OCompIneq}), so they match with the $\eta\neq k$ case. They also shift identically since the regions are identical to the $\eta \neq k$ case.

\section{Singular Solution in the Original Variables}

Now that the cases have been generally solved for, we ensure consistency with the original balance equations.  

\subsection{Case $\eta \neq k$}
When we return to the original variables, the delta-shock solution is represented in the following way
\begin{equation}
    \label{geneq_2}
    \begin{aligned}
    \big(\rho,u\big)\big(x,t\big) =
    \begin{cases}
        \bigg(\rho_Le^{kt},\bigg(u_L+\frac{\beta}{\eta-k}\bigg)e^{(\eta-k)t}-\frac{\beta}{\eta-k}\bigg) \text{,     } & x<x\big(t\big)
        \\ \\ 
        \bigg(\bar{\omega}\big(t\big)\delta\big(x-x\big(t\big)\big), u_{\delta}\big(t\big)\bigg) \text{,  } & x=x\big(t\big)
        \\ \\ 
        \bigg(\rho_Re^{kt},\bigg(u_R+\frac{\beta}{\eta-k}\bigg)e^{(\eta-k)t}-\frac{\beta}{\eta-k}\bigg) \text{,    } & x>x\big(t\big), \text{ }
    \end{cases}
    \end{aligned}
\end{equation}
where $\bar{\omega}(t)=\omega_1(t)e^{kt}.$
Converting (\ref{shock}), (\ref{n=kRanHugRel_1}),(\ref{n=kRanHugRel_2}):
\begin{equation}
\label{RankineP1}
\begin{aligned}
\begin{split}
    \frac{d\omega_1}{dt}&=e^{-kt}\frac{d\Bar{\omega}}{dt}-\bar{\omega}ke^{-kt}\\&=-\big(\rho_L-\rho_R\big)u_{\delta}\big(t\big)+\rho_L\bigg(u_L+\frac{\beta}{\eta-k}\bigg)e^{(\eta-k)t}-\frac{\rho_L\beta}{\eta-k}-A\big(\rho_Le^{kt}\big)^{\gamma+1}e^{(\eta-k)t}\\
    &-\rho_R\bigg(u_R+\frac{\beta}{\eta-k}\bigg)e^{(\eta-k)t}+\frac{\rho_R\beta}{\eta-k}-A\big(\rho_Re^{kt}\big)^{\gamma+1}e^{(\eta-k)t},
\end{split}
\end{aligned}
\end{equation} 
\begin{equation}
\label{RankineP2}
\begin{aligned}
\begin{split}
    \frac{d}{dt}\bigg(\omega_1\bigg(\omega_\delta+\frac{\beta}{\eta-k}\bigg)\bigg)&=\frac{d}{dt}\bigg(\bar{\omega}e^{-kt}\bigg(u_\delta+\frac{\beta}{\eta-k}\bigg)e^{-(\eta-k)t}\bigg)\\&=-\bigg(\rho_L\bigg(u_L+\frac{\beta}{\eta-k}\bigg)-\rho_{R}\bigg(u_R+\frac{\beta}{\eta-k}\bigg)\bigg)u_{\delta}(t)\\
    &+\bigg(\bigg(u_L+\frac{\beta}{\eta-k}\bigg)^2\rho_Le^{(\eta-k)t}-\bigg(u_L+\frac{\beta}{\eta-k}\bigg)e^{(\eta-k)t}A\big(v_Le^{kt}\big)^{\gamma+1}\\
    &-\frac{\beta}{\eta-k}\bigg(u_L+\frac{\beta}{\eta-k}\bigg)\rho_L+\frac{\beta}{\eta-k}\bigg(u_R+\frac{\beta}{\eta-k}\bigg)\rho_R\\
    &-\bigg(u_R+\frac{\beta}{\eta-k}\bigg)^2\rho_Re^{(\eta-k)t}+\bigg(u_R+\frac{\beta}{\eta-k}\bigg)e^{(\eta-k)t}A\big(v_Re^{kt}\big)^{\gamma+1}\bigg),
\end{split}
\end{aligned}
\end{equation}
produces
\begin{equation}\label{RankineO}
\frac{dx}{dt}=u_\delta,
\end{equation}
\begin{equation}
    \label{omegabar}
    \frac{d\bar{\omega}}{dt}=k\bar{\omega}\big(t\big)-\big[\rho\big]_{\text{jump}}u_{\delta}\big(t\big)+\big[\rho u\big]_{\text{jump}}-A\big[\rho^{\gamma+1}\big]_{\text{jump}}e^{nt},
\end{equation}
and 
\begin{equation}
\label{Rankine2}
    \frac{d}{dt}\bigg(\bar \omega u_\delta\bigg)=\eta\bar{\omega}\big(t\big)u_\delta+\beta\bar{\omega}\big(t\big)-\big[\rho u\big]_{\text{jump}}u_\delta+\big[\rho u^2\big]_{\text{jump}}-A\big[\rho^{\gamma+1}u\big]_{\text{jump}}e^{\eta t}.
\end{equation}
Substituting $u$ and $\rho$ into the delta-shock definition to verify that the solution satisfies the equations in the sense of distributions yields
\begin{equation}
\label{RankineDef}
\begin{cases}
    \bigg\langle \rho,\phi_t\bigg\rangle+\bigg\langle\rho u-A\rho^{\gamma+1}e^{\eta t},\phi_x\bigg\rangle=-\bigg\langle k\rho,\phi\bigg\rangle\\\bigg\langle \rho u,\phi_t\bigg\rangle+\bigg\langle\rho u^2-Au\rho^{\gamma+1}e^{\eta t},\phi_x\bigg\rangle=-\bigg\langle \eta\rho u,\phi\bigg\rangle-\bigg\langle \beta\rho,\phi\bigg\rangle.
\end{cases}
\end{equation}
Only the proof of the second equality is presented. The first can be shown to hold by a similar argument.  
Let 
\begin{equation}
\label{RankineQR}
\begin{cases}
    Q=\bigg\langle \rho u,\phi_t\bigg\rangle+\bigg\langle\rho u^2-Au\rho^{\gamma+1}e^{\eta t},\phi_x\bigg\rangle\\R=-\bigg\langle \eta\rho u,\phi\bigg\rangle-\bigg\langle \beta\rho,\phi\bigg\rangle
\end{cases}
\end{equation}
then
\begin{equation}
\label{RankineQ}
\begin{split}
    Q=&\int_0^\infty \int_{-\infty}^{x(t)}\rho_Le^{kt}\bigg(\bigg(u_L+\frac{\beta}{\eta-k}\bigg)e^{(\eta-k)t}-\frac{\beta}{\eta-k}\bigg)\phi_tdxdt\\&+\int_0^{\infty}\int_{x(t)}^\infty\rho_Re^{kt}\bigg(\bigg(u_R+\frac{\beta}{\eta-k}\bigg)e^{(\eta-k)t}-\frac{\beta}{\eta-k}\bigg)\phi_tdxdt\\
    &+\int_0^\infty \int_{-\infty}^{x(t)}\rho_Le^{kt}\bigg(\bigg(u_L+\frac{\beta}{\eta-k}\bigg)e^{(\eta-k)t}-\frac{\beta}{\eta-k}\bigg)^2\phi_xdxdt\\&+\int_0^\infty\int_{x(t)}^\infty \rho_Re^{kt}\bigg(\bigg(u_R+\frac{\beta}{\eta-k}\bigg)e^{(\eta-k)t}-\frac{\beta}{\eta-k}\bigg)^2\phi_xdxdt\\&-A\int_0^{\infty}\int_{-\infty}^{x(t)}\big(\rho_Le^{kt}\big)^{\gamma+1}\bigg(\bigg(u_L+\frac{\beta}{\eta-k}e^{(\eta-k)t}\bigg)-\frac{\beta}{\eta-k}\bigg)e^{\eta t}\phi_xdxdt\\&-A\int_0^\infty\int_{-\infty}^{x(t)}\big(\rho_Re^{kt}\big)^{\gamma+1}\bigg(\bigg(u_R+\frac{\beta}{\eta-k}\bigg)e^{(\eta-k)t}-\frac{\beta}{\eta-k}\bigg)e^{\eta t}\phi_xdxdt\\&+\int_0^\infty \bar{\omega} u_\delta\bigg(\phi_t+u_\delta\phi_x\bigg)dt.
\end{split}
\end{equation}
If we assume that $\frac{dx}{dt}>0$ for $t\in\mathbb{R}^+$ (we employ a similar argument in the case $\frac{dx}{dt}<0$), then an inverse of $x(t)$ exists.  Thus, 
\begin{equation}
\label{RankineQ1}
\begin{split}
    Q=&\int_0^\infty\int_{t(x)}^\infty \rho_Le^{kt}\bigg(\bigg(u_L+\frac{\beta}{\eta-k}\bigg)e^{(\eta-k)t}-\frac{\beta}{\eta-k}\bigg)\phi_tdtdx\\&+\int_0^\infty \int_0^{t(x)} \rho_Re^{kt}\bigg(\bigg(u_R+\frac{\beta}{\eta-k}\bigg)e^{(\eta-k)t}-\frac{\beta}{\eta-k}\bigg)\phi_tdtdx\\&+\int_0^\infty \int_{-\infty}^{x(t)}\rho_Le^{kt}\bigg(\bigg(u_L+\frac{\beta}{\eta-k}\bigg)e^{(\eta-k)t}-\frac{\beta}{\eta-k}\bigg)^2\phi_xdxdt\\&+\int_0^\infty\int_{x(t)}^\infty \rho_Re^{kt}\bigg(\bigg(u_R+\frac{\beta}{\eta-k}\bigg)e^{(\eta-k)t}-\frac{\beta}{\eta-k}\bigg)^2\phi_xdxdt\\&-A\int_0^{\infty}\int_{-\infty}^{x(t)}\big(\rho_Le^{kt}\big)^{\gamma+1}\bigg(\bigg(u_L+\frac{\beta}{\eta-k}e^{(\eta-k)t}\bigg)-\frac{\beta}{\eta-k}\bigg)e^{\eta t}\phi_xdxdt\\&-A\int_0^\infty\int_{x(t)}^{\infty}\big(\rho_Re^{kt}\big)^{\gamma+1}\bigg(\bigg(u_R+\frac{\beta}{\eta-k}\bigg)e^{(\eta-k)t}-\frac{\beta}{\eta-k}\bigg)e^{\eta t}\phi_xdxdt\\&+\int_0^\infty \bar{\omega} u_\delta d\phi.
\end{split}
\end{equation}
After an integration by parts, we deduce
\begin{equation}
\label{RankineQ1Simp}
\begin{split}
    Q=&-\eta \int_0^\infty \int_{t(x)}^\infty \rho_Le^{\eta t}\bigg(\bigg(u_L+\frac{\beta}{\eta-k}\bigg)e^{(\eta-k)t}-\frac{\beta}{\eta-k}\bigg)\phi dtdx\\&-\beta\int_0^\infty\int_{t(x)}^\infty\rho_Le^{kt}\phi dtdx\\&-\beta\int_0^\infty\int_{0}^{t(x)}\rho_Re^{kt}\phi dtdx\\&-\eta \int_0^\infty\int_0^{t(x)}\rho_Re^{kt}\bigg(\bigg(u_R+\frac{\beta}{\eta-k}\bigg)e^{(\eta-k)t}-\frac{\beta}{\eta-k}\bigg)\phi dtdx\\&+\int_0^\infty N\big(t\big)\phi\big(x\big(t\big),t\big)dt,
\end{split}
\end{equation}
where
\begin{equation}
\label{RankineN}
\begin{split}
    N(t)=&\bigg(\rho_+e^{kt}\bigg(\bigg(u_++\frac{\beta}{\eta-k}\bigg)e^{(\eta-k)t}-\frac{\beta}{\eta-k}\bigg)-\rho_-e^{kt}\bigg(\bigg(u_-+\frac{\beta}{\eta-k}\bigg)e^{(\eta-k)t}-\frac{\beta}{\eta-k}\bigg)\bigg)u_\delta\\&+\rho_-e^{kt}\bigg(\bigg(u_-+\frac{\beta}{\eta-k}\bigg)e^{(\eta-k)t}-\frac{\beta}{\eta-k}\bigg)^2-Ae^{\eta t}\big(\rho_-e^{kt}\big)^{\gamma+1}\bigg(\bigg(u_-+\frac{\beta}{\eta-k}\bigg)e^{(\eta-k)t}-\frac{\beta}{\eta-k}\bigg)\\&-\rho_+e^{kt}\bigg(\bigg(u_++\frac{\beta}{\eta-k}\bigg)e^{(\eta-k)t}-\frac{\beta}{\eta-k}\bigg)^2+Ae^{\eta t}\big(\rho_+e^{kt}\big)^{\gamma+1}\bigg(\bigg(u_++\frac{\beta}{\eta-k}\bigg)e^{(\eta-k)t}-\frac{\beta}{\eta-k}\bigg)\\
    &-\frac{d}{dt}\bigg(\bar{\omega}\big(t\big)u_{\delta}\bigg)\\
    =& -\big[\rho u\big]_{\text{jump}}u_{\delta}+\bigg[\rho u^2-Ae^{\eta t}\rho^{\gamma+1}u\bigg]_{\text{jump}}-\frac{d}{dt}\bigg(\bar{\omega}\big(t\big)u_{\delta}\bigg)\\
    \myeq &-\eta \bar{\omega}u_\delta-\beta \bar{\omega}\big(t\big).
\end{split}
\end{equation}
We can easily conclude that $Q=R$ as desired, so the equations are satisfied in the sense of distributions. Now consider the following case:
\begin{figure}[H]
\begin{center}
\begin{tikzpicture}
  \begin{axis}[
      axis x line=middle,
      axis y line=middle,
      ymax=2.5,
      ymin=-1,
      xmax=2,
      xmin=-1,
      domain=0:10,
      xtick={0, 1, ..., 5},
      ytick={-5, -4, ..., 5},
      samples=1001,
      yticklabels={, , , , , , $t_1$,,},
      xticklabels={0,$x_0$,},
      xlabel = \(x\),
      ylabel = {\(t\)}
    ]
    \addplot [
    color=black,
    ] {-((1-x^2)^(1/2))+1};
    \addplot [
    color=black,
    ] {((1-x^2)^(1/2))+1};
    \addplot [
    domain=-1:0,
    color=black, very thick
    ] {2};
    \addplot [
    domain = .866:1.3,
    color=black,
    ] {1.5};
    \addplot [
    domain = .866:1.3,
    color=black,
    ] {.5};
    \draw (250,250) node[draw] {$ \  \ \ t = t_2(x)$};
    \node at (axis cs:-0.13,2.2)(B){\color{black}$t_0$};
    \draw (250,150) node[draw] {$\ \ \ t = t_1(x)$};
    \addplot [
    color = black, very thick, dotted,
    domain = 0:1,
    ] {1};
  \end{axis}
\end{tikzpicture}
\end{center}
\caption{\textit{Graph of an example where $x'(t_1)=0$}}
    \label{fig:WellDefinedInverse}
\end{figure}
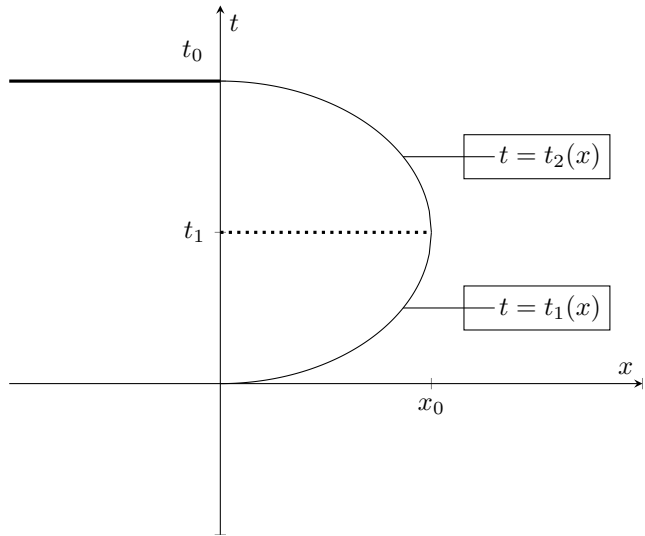
If $x=x(t)$ is a curve, as shown above, the earlier proof needs to be modified. If there are more points at which $x'(t)=0,$ the proof is similar. The proof can be modified by breaking up $Q$ as follows

\begin{equation}
\label{RankineQ2}
\begin{split}
    Q=&\int_{-\infty}^0\int_0^{t(x)}\rho_Le^{kt}\bigg(\bigg(u_L+\frac{\beta}{\eta-k}\bigg)e^{(\eta-k)t}-\frac{\beta}{\eta-k}\bigg)\phi_tdtdx\\&+\int_{-\infty}^0\int_{t_2(x)}^\infty \rho_Re^{kt}\bigg(\bigg(u_R+\frac{\beta}{\eta-k}\bigg)e^{(\eta-k)t}-\frac{\beta}{\eta-k}\bigg)\phi_tdtdx\\&+\int_0^{x_0}\int_0^{t_1(x)}\rho_Re^{kt}\bigg(\bigg(u_R+\frac{\beta}{\eta-k}\bigg)e^{(\eta-k)t}-\frac{\beta}{\eta-k}\bigg)\phi_tdtdx\\&+\int_0^{x_0}\int_{t_1(x)}^{t_2(x)}\rho_Le^{kt}\bigg(\bigg(u_L+\frac{\beta}{\eta-k}\bigg)e^{(\eta-k)t}-\frac{\beta}{\eta-k}\bigg)\phi_tdtdx\\&+\int_0^{x_0}\int_{t_2(x)}^{\infty}\rho_Re^{kt}\bigg(\bigg(u_R+\frac{\beta}{\eta-k}\bigg)e^{(\eta-k)t}-\frac{\beta}{\eta-k}\bigg)\phi_tdtdx\\&+\int_{x_0}^{\infty}\int_{0}^{\infty}\rho_Re^{kt}\bigg(\bigg(u_R+\frac{\beta}{\eta-k}\bigg)e^{(\eta-k)t}-\frac{\beta}{\eta-k}\bigg)\phi_tdtdx\\&+\int_0^{\infty}\int_{-\infty}^{x(t)}\rho_Le^{kt}\bigg(\bigg(u_L+\frac{\beta}{\eta-k}\bigg)e^{(\eta-k)t}-\frac{\beta}{\eta-k}\bigg)^2\phi_xdxdt\\&+\int_{0}^{\infty}\int_{x(t)}^{\infty}\rho_Re^{kt}\bigg(\bigg(u_R+\frac{\beta}{\eta-k}\bigg)e^{(\eta-k)t}-\frac{\beta}{\eta-k}\bigg)^2\phi_xdxdt\\&-A\int_0^{\infty}\int_{-\infty}^{x(t)}\big(\rho_Le^{kt}\big)^{\gamma+1}\bigg(\bigg(u_L+\frac{\beta}{\eta-k}\bigg)e^{(\eta-k)t}-\frac{\beta}{\eta-k}\bigg)^2\phi_xdxdt\\&-A\int_{0}^{\infty}\int_{x(t)}^{\infty}\big(\rho_Re^{kt}\big)^{\gamma+1}\bigg(\bigg(u_R+\frac{\beta}{\eta-k}\bigg)e^{(\eta-k)t}-\frac{\beta}{\eta-k}\bigg)^2\phi_xdxdt\\&+\int_0^\infty \bar{\omega}u_\delta d\phi
\end{split}
\end{equation}
After an integration by parts and change of variables by using (\ref{RankineO}), we once again conclude that $Q=R.$
For a strictly overcompressive delta-shock solution connecting a left state $(\rho_L,u_L)$ and a right state $(\rho_R,u_R)$ we require 

\begin{equation}
\label{OCEQ}
    u_R-A\rho_R^{\gamma}e^{\eta t}<u_\delta<u_L-A\rho_L^{\gamma}e^{\eta t}\big(\gamma+1\big).
\end{equation}



\subsection{Case $\eta =k$}
In this case, the delta-shock solution is represented in the following way
\begin{equation}
    \label{geneq}
    \begin{aligned}
    \big(\rho,u\big)\big(x,t\big) =
    \begin{cases}
        \bigg(\rho_Le^{kt},u_L+\beta t\bigg) \text{,     } & x<x\big(t\big)
        \\ \\ 
        \bigg(\bar{\omega}\big(t\big)\delta\big(x-x\big(t\big)\big), u_{\delta}\big(t\big)\bigg) \text{,  } & x=x\big(t\big)
        \\ \\ 
        \bigg(\rho_Re^{kt},u_R+\beta t\bigg) \text{,    } & x>x\big(t\big), \text{ }
    \end{cases}
    \end{aligned}
\end{equation}
where $\bar{\omega}(t)=\omega_1(t)e^{kt}.$
From (\ref{shock_2}), (\ref{n=kRanHugRel}), we get
\begin{equation}\label{Rankine_3}
\frac{dx}{dt}=\sigma\big(t\big)=u_\delta,
\end{equation}
and
\begin{equation}
\label{RankineP_1}
\begin{cases}
\begin{aligned}
    \frac{d\omega_1}{dt}&=e^{-kt}\frac{d\Bar{\omega}}{dt}-\bar{\omega}ke^{-kt}\\&=-\big(v_L-v_R\big)\sigma\big(t\big)+\rho_L\big(w_L+\beta t\big)-A\big(\rho_Le^{kt}\big)^{\gamma+1}-\rho_R\big(w_R+\beta t\big)+A\big(\rho_Re^{kt}\big)^{\gamma+1}\\
    &\Rightarrow \frac{d\bar{\omega}}{dt}=k\bar{\omega}\big(t\big)-\big[\rho\big]_{\text{jump}}+\bigg[\rho u-A\rho^{\gamma+1}e^{kt}\bigg]_{\text{jump}},
    \end{aligned}
\\ \\
\begin{aligned}
\frac{d}{dt}\bigg(\bar{\omega}e^{-kt}\big(u_{\delta}-\beta t\big)\bigg)&=-\big(v_Lw_L-v_R w_R\big)u_{\delta}\\
&+v_L\big(w_L+\beta t\big)w_L-Av_L^{\gamma+1}e^{kt(\gamma+1)} w_L-v_R\big(w_R+\beta t\big)w_R+Av_R^{\gamma+1}e^{kt(\gamma+1)} w_R\\
&\Rightarrow \frac{d}{dt}\big(\bar{\omega}u_{\delta}\big)=k\bar{\omega}\big(t\big)u_{\delta}-\big[\rho u\big]_{\text{jump}}u_{\delta}+\bigg[\rho u^2-A\rho^{\gamma+1}e^{kt}u\bigg]_{\text{jump}}+\bar{\omega}\beta,
\end{aligned}
\end{cases}
\end{equation} 
respectively. All remaining proofs for $\eta = k$ follow precisely the $\eta \neq k$ case, including those for the definitions of delta-shocks and overcompressive regions.

\section{Conclusion}
In this work, we studied the Riemann problem of a non-symmetric Keyfitz-Kranzer type system with varying generalized Chaplygin gas. While there is a substantial body of literature on the Chaplygin gas and, more recently, on the varying Chaplygin gas, see Li \cite{Li}, our study takes a unique approach and combines various models. We address the open question of whether classical and non-classical (delta-shocks) solutions are possible in the presence of a power $\gamma$ in the density.

\vspace{5mm}

We provide an affirmative answer by deriving various regions in four cases (depending on the sign of $k(\gamma+1)$ and whether $\gamma$ is less or greater than $-1$), where the Riemann problem can be solved classically (by using a 1-shock, 2-rarefaction, 2-contact discontinuity) or non-classically (by using a combination of classical waves and a delta-shock, or solely a delta-shock). We observe that these regions shift in time. Therefore, a Riemann problem with a given left and right state can have different solutions over several time intervals. We also prove that the singular solution (which involves a delta-shock) satisfies our system in the sense of distributions.  More generally, the results highlight the challenge of solving the Riemann problem for a non-autonomous system of balance laws (presence of source terms) due to the lack of self-similarity and direct dependence on time of the wave curves, which causes region shifts.

\vspace{5mm}

Lastly, our robust numerical evidence indicates the existence of regions where the solution consists of a 2-rarefaction followed by a 2-contact discontinuity (in that specific order), which we have not verified analytically, and regions where the solutions consist of a combination of a classical wave and a delta-shock. We verified the feasibility of the Local Lax-Friedrichs method for time-dependent solutions by manipulating key parameters to study changes in time. Future work will pursue these and other questions, such as how these Riemann solutions can be used as building blocks in solving general Cauchy problems.

\vspace{5mm}

{\bf Acknowledgment.} This work is supported by the National Science Foundation under Grant Number DMS-2349040 (PI: Tsikkou). Any opinions, findings, and conclusions or recommendations expressed in this
material are those of the authors and do not necessarily
reflect the views of the National Science Foundation. \\

\noindent
The authors thank Barbara Lee Keyfitz for providing MATLAB code that served as a basis for the numerical analysis done in this paper. \\

\noindent
{\bf Data Availability.} The data that support the findings of
this study are available from the corresponding author
upon reasonable request.

\end{document}